\documentclass[11pt,a4paper]{amsart}
\usepackage[utf8]{inputenc}
\usepackage{amsmath}
\usepackage{amsthm}
\usepackage{amsfonts}
\usepackage{amssymb}
\usepackage[english]{babel}
\usepackage{MnSymbol}
\usepackage{cancel}
\usepackage{mathtools}
\usepackage{comment}
\usepackage{paralist}
\usepackage{caption}
\usepackage{enumerate}

\usepackage{lipsum}
\newcommand\blfootnote[1]{%
  \begingroup
  \renewcommand\thefootnote{}\footnote{#1}%
  \addtocounter{footnote}{-1}%
  \endgroup}

\newtheorem{The}{Theorem}[section]
\newtheorem*{The*}{Theorem}
\newtheorem{thmx}{Theorem}

\newtheorem{Lem}[The]{Lemma}
\newtheorem{Def}[The]{Definition}
\newtheorem{Prop}[The]{Proposition}
\newtheorem{Example}[The]{Example}
\newtheorem{Remark}[The]{Remark}
\newtheorem{Cor}[The]{Corollary}
\newtheorem*{claim}{Claim}

\numberwithin{equation}{section}

\newcommand{\Ht}{\mathbb{H}^2}
\newcommand{\g}{\mathfrak{g}}
\newcommand{\p}{\mathfrak{p}}
\newcommand{\kk}{\mathfrak{k}}

\newcommand{\am}{\mathfrak{a}}
\newcommand{\apo}{\am^{+}_1}
\newcommand{\ap}{\am^{+}}
\newcommand{\at}{\am_{\tm}}

\newcommand{\op}{\texttt{\hspace{1mm}op\hspace{1mm}}}
\newcommand{\nop}{\texttt{\hspace{1mm}\cancel{op}\hspace{1mm}}}

\newcommand{\sm}{\sigma}
\newcommand{\tm}{\tau}

\newcommand{\Ftm}{\mathrm{Flag}_{\tau}}
\newcommand{\Ftim}{\mathrm{Flag}_{\iota\tau}}
\newcommand{\Ftsm}{\mathrm{Flag}_{\sigma}}
\newcommand{\Ftmo}{\mathrm{Flag}_{\tau_1}}
\newcommand{\Ftmt}{\mathrm{Flag}_{\tau_2}}

\newcommand{\ty}{\mathrm{typ}}

\newcommand{\id}{\mathrm{id}}

\newcommand{\Af}{\mathcal{A}_{\xi}}

\newcommand{\At}{\mathcal{A}_{\tm}}
\newcommand{\Atop}{\mathcal{A}^{\texttt{op}}_{\tm}}
\newcommand{\Asop}{\mathcal{A}^{\texttt{op}}_{\sm}}

\newcommand{\Ato}{\mathcal{A}_{\tm_1}}
\newcommand{\Att}{\mathcal{A}_{\tm_2}}

\newcommand{\rx}{\rho_{x,y}}
\newcommand{\ry}{\rho_{y,x}}

\newcommand{\st}{\mathrm{stab}}
\newcommand{\inn}{\mathrm{int}}
\newcommand{\Sc}{\textsc{S}}

\newcommand{\cro}{\mathrm{cr}}
\newcommand{\Ima}{\mathrm{Im}}
\newcommand{\CAT}{\mathrm{CAT}(0)}

\newcommand{\R}{\mathbb{R}}
\newcommand{\N}{\mathbb{N}}

\newcommand{\di}{\partial_{\infty}}
\newcommand{\Di}{\Delta_{\infty}}

\newcommand{\Res}{\mathrm{Res}}
\newcommand{\rk}{\mathrm{rk}}

\newcommand{\bl}{\begin{Lem}}
\newcommand{\el}{\end{Lem}}

\newcommand{\bd}{\begin{Def}}
\newcommand{\ed}{\end{Def}}

\newcommand{\bt}{\begin{The}}
\newcommand{\et}{\end{The}}

\newcommand{\bp}{\begin{Prop}}
\newcommand{\ep}{\end{Prop}}

\newcommand{\bc}{\begin{Cor}}
\newcommand{\ec}{\end{Cor}}

\begin{document}

\title[Cross ratios on symmetric spaces and Euclidean buildings]{Cross ratios on boundaries of symmetric spaces and Euclidean buildings}
\author[Jonas Beyrer]{Jonas Beyrer}
\date{}

\maketitle

\textbf{Abstract:} We generalize the natural cross ratio on the ideal boundary of a rank one symmetric spaces, or even $\mathrm{CAT}(-1)$ space, to higher rank symmetric spaces and (non-locally compact) Euclidean buildings. We obtain vector valued cross ratios defined on simplices of the building at infinity. We show several properties of those cross ratios; for example that (under some restrictions) periods of hyperbolic isometries give back the translation vector. In addition, we show that cross ratio preserving maps on the chamber set are induced by isometries and vice versa - motivating that the cross ratios bring the geometry of the symmetric space/Euclidean building to the boundary.
\blfootnote{The author was support by the SNF grant \text{200020{\_}175567}}

\section{Introduction}

Cross ratios on boundaries are a crucial tool in hyperbolic geometry and more general negatively curved spaces. In this paper we show that we can generalize these cross ratios to (the non-positively curved) symmetric spaces of higher rank and thick Euclidean buildings with many of the properties of the cross ratio still valid.

On the boundary $\di \Ht$ of the hyperbolic plane $\Ht$ there is naturally a multiplicative cross ratio defined by $\cro_{\Ht}(z_1,z_2,z_3,z_4)=\frac{z_1-z_2}{z_1-z_4}\frac{z_3-z_4}{z_3-z_2}$ when considering $\Ht$ in the upper half space model, i.e. $\di \Ht=\R\cup\{\infty\}$. This cross ratio plays an essential role in hyperbolic geometry. For example it characterizes the isometry group by its boundary action and therefore allows to study the geometry of the space from its boundary; which is an important perspective in hyperbolic geometry.

This cross ratio can be generalized in a way broader context, namely $\mathrm{CAT}(-1)$ spaces \cite{Bourdon-EM}: Let $\di Y$ be the ideal boundary of a $\mathrm{CAT}(-1)$ space $Y$, $x,y\in \di Y$ and $o\in Y$. Then the \emph{Gromov product} $(\cdot |\cdot)_o:\di Y^2\to [0,\infty]$ is defined by $(x|y)_o=\lim_{t\to\infty}t- \frac{1}{2}d(\gamma_{ox}(t),\gamma_{oy}(t))$, where $\gamma_{ox},\gamma_{oy}$ are the unique unit speed geodesics from $o$ to $x,y$, respectively. 
Then a \emph{additive cross ratio} $\cro_{\di Y}: \mathcal{A}\subset\di Y^4\to [0,\infty]$ is defined by $\cro_{\di Y}(x,y,z,w):=-(x|y)_o-(z|w)_o+(x|w)_o+(z|y)_o$ for all $(x,y,z,w)\in \di Y^4$ with no entry occurring three or four times;\footnote{This would correspond to considering $\log |\cro_{\Ht}|$ for the hyperbolic plane}
which is independent of the basepoint. By construction $\cro_{\di Y}$ has several symmetries with respect to $(\R,+)$.
In analogy to the hyperbolic plane, maps $f:\di Y\to\di Y$ that leave $\cro_{\di Y}$ under the diagonal action invariant are called \emph{Moebius maps}. It follows from the definition of the cross ratio together with the basepoint independence that isometries are Moebius maps when restricted to the boundary.

The cross ratios $\cro_{\di Y}$ and Moebius maps have been proven to be very useful in hyperbolic geometry. For example Bourdon \cite{Bourdon} has shown that Moebus maps of rank one symmetric spaces extend uniquely to isometric embeddings of the interior, and with this he gave a new proof of Hamenst\"adt's 'entropy against curvature' theorem \cite{Hamenstaedt-Ann}. 
Otal \cite{Otal-Ann} has (implicitly) shown that Moebius bijections on boundaries of universal covers of closed negatively-curved surfaces can be uniquely extended to isometries; which yields that marked length spectrum rigidity holds for those manifolds - a prominent conjecture formulated in \cite{Burns-Katok}. See \cite{Dalbo-Kim,Kim-Topology01,Kim-Topology04} for more results in that context. 
Moreover, there is a close relation between the cross ratio on the boundary of the universal cover of a closed negatively curved manifold and the quasi-conformal structure on the boundary, and to dynamical properties of the
geodesic flow; see e.g. \cite{Ledrappier}.

On the boundary $\di \tilde{S}$ of the universal cover of a closed surface $S$ there are many other cross ratios, besides the above constructed one, that parametrize classical objects associated to the surface; such as simple closed curves, measured laminations, points of Teichm\"uller space \cite{Bonahon-IM}, Hitchin representations \cite{Labourie-IHES} and positively ratioed representations \cite{Martone-Zhang}\footnote{We will see that the cross ratios associated to Hitchin representations and positively ratioed representations arise as pullbacks (under the natural boundary map) of cross ratios that we construct in this paper.} - to name a few.\medskip

This prominence and importance of cross ratios in negative curvature motivates us to ask if such objects also exists for non-positively curved spaces and how much information about the geometry they carry. 

There is already some work done in this context. In \cite{Charney-Murray} a coarse cross ratio for arbitrary $\CAT$ spaces on some subset of the boundary has been constructed. In \cite{Beyrer-Fioravanti-Incerti} there is a cross ratio defined on the Roller boundary of a $\CAT$ cube complexes, using essentially the combinatorial structure of the space. In those works Moebius (respectively quasi-Moebius) bijections are connected to isometries (respectively quasi-isometries).

In this paper we will construct cross ratios for symmetric spaces and Euclidean buildings, which will generalize the cross ratios of $\mathrm{CAT}(-1)$ spaces. 
There is little need to explain the importance of symmetric spaces in differential geometry and related areas. However, we want to point out that the study of symmetric spaces has recently gained renewed prominence in the active field of research of Anosov representations and Anosov subgroups (e.g. \cite{Labourie-IM}, \cite{KLP}, \cite{Guichard-Wienhard-IM} and many more). We will see that the cross ratios we construct are connected to the study of those (e.g. \cite{Labourie-IHES}, \cite{Martone-Zhang}) and hence we hope for applications of our work in this area.

Euclidean buildings arise in many different areas of mathematics. See \cite{Ji} for an overview of some applications. Probably most prominently they arise in the study of algebraic groups and geometric group theory; they have also been a crucial tool in the proof of quasi-isometric rigidity of symmetric spaces \cite{Kleiner-Leeb} (extending Mostow-Prasad rigidity) - to name a few.\medskip
 
We will denote by $M$ either a symmetric space or a thick Euclidean building. It is well known that the ideal boundary $\di M$ has naturally the structure of a spherical building $\Di M$. Therefore there is a \emph{type map} $\ty:\di M\to \sm$ with $\sm$ the closed fundamental chamber of the spherical Coxeter complex associated to $M$.
Then we show that to each \emph{type} $\xi\in\sm$ there is $\iota\xi\in\sm$ such that the Gromov product (defined exactly as for $\mathrm{CAT}(-1)$ spaces) restricted to the set $\ty^{-1}(\xi)\times \ty^{-1}(\iota\xi)$ is generically finite. Thus we get a generically defined \emph{additive} cross ratio on $(\ty^{-1}(\xi)\times \ty^{-1}(\iota\xi))^2$ in the same way as for $\mathrm{CAT}(-1)$ spaces. We can show that this cross ratio is independent of the choice of basepoint; and denote it by $\cro_{\xi}$.

Let $\tm$ be a face of the simplex $\sm$, $\inn(\tm)$ the interior of $\tm$ and $\xi\in\inn(\tm)$. Moreover, we denote by $\Ftm(M)\subset\Di M$ the set of simplices of the building at infinity of type $\tm$ (i.e. those simplices that are mapped to $\tm$ under $\ty$); in particular $\Ftsm(M)$ is the chamber set of the building at infinity. Then one can naturally identify $\ty^{-1}(\xi)$ with $\Ftm(M)$ and in the same way $\ty^{-1}(\iota\xi)$ with $\Ftim(M)$. 
Therefore we immediately get a cross ratio $\cro_{\xi}:\mathcal{A}_{\tm}\subset (\Ftm(M)\times \Ftim(M))^2\to [-\infty,\infty]$, which by construction has similar symmetries as the additive one on $\mathrm{CAT}(-1)$ spaces - for $\mathcal{A}_{\tm}$ see equation \eqref{eq def domain of definition}, for the symmetries see equation \eqref{eq symmetries of cross ratios}.

Clearly, we get a whole collection of cross ratios defined on the set $\mathcal{A}_{\tm}$ which is parametrized by $\xi\in\inn(\tm)$. Then we show that we can put together this collection to a single vector valued cross ratio $\cro_{\tm}$ with the same symmetries, and values in the Coxeter complex associated to $M$. We will see that the vector valued cross ratio is the natural object to consider; we can connect the so called \emph{period} $\cro_{\sm}(g^{-},g\cdot x,g^{+},x)$ of a hyperbolic element $g\in Iso(M)$ (with attractive and repulsive fixed points $g^{\pm}\in\Ftsm(X)$ and generic $x\in\Ftsm(X)$) to the translation vector of $g$ along the unique maximal flat joining $g^{-}$ and $g^{+}$, and we give a `nice' geometric interpretation of the vector valued cross ratio.

Let $M_1,M_2$ be either two symmetric spaces or two thick Euclidean buildings. Let $\sm_1,\sm_2$ be the according fundamental chambers of the spherical Coxeter complexes and let $\xi_i\in\inn(\sm_i)$ be two types. Let $f:\Ftsm(M_1)\to\Ftsm(M_2)$ be surjective. The map $f$ is called \emph{$\xi_1$-Moebius bijection}, if $\cro_{\xi_1}(x,y,z,w)=\cro_{\xi_2}(f(x),f(y),f(z),f(w))$ for all $(x,y,z,w)\in\mathcal{A}_{\sm_1}$, and \emph{$\sm_1$-Moebius bijection}, if $\cro_{\sm_1}(x,y,z,w)=\cro_{\sm_2}(f(x),f(y),f(z),f(w))$ for all $(x,y,z,w)\in\mathcal{A}_{\sm_1}$. Moreover, we call a locally compact Euclidean building with discrete translation group a  \emph{combinatorial Euclidean building} and a Euclidean building \emph{thick} if and only if the building at infinity is thick. Then we show the following:

\begin{thmx}
Let $M_1,M_2$ be either symmetric spaces or thick combinatorial Euclidean buildings and $\xi_1\in\inn(\sm_1)$. If $M_1,M_2$ are irreducible, then every $\xi_1$-Moebius bijection $f:\Ftsm(M_1)\to\Ftsm(M_2)$ can be extended to an isometry $F:M_1\to M_2$. If none of the spaces is a Euclidean cone over a spherical building, then this extension is unique. If $M_1,M_2$ are reducible one can rescale the metric of $M_1$ on irreducible factors - denote this space by $\hat{M}_1$ - such that $f$ can be extended to an isometry $F:\hat{M}_1\to M_2$.
\end{thmx}

\begin{thmx}\label{thm B} Let $E_1,E_2$ be thick (non-locally compact) Euclidean buildings. Then for every $\sm_1$-Moebius bijection $f:\Ftsm(E_1)\to\Ftsm(E_2)$ one can rescale the metric of $E_1$ on irreducible factors - denote this space by $\hat{E}_1$ - such that $f$ can be extended to an isometry $F:\hat{E}_1\to E_2$. If none of the irreducible factors is a Euclidean cone over a spherical building, then $f$ can be extended to an isometry $F:E_1\to E_2$ (without rescaling the metric).
\end{thmx}

We remark that essentially by definition of the cross ratio every isometry gives rise to a Moebius bijection. Therefore these theorems show that the cross ratios - at least for the chamber set of the building at infinity - carry a lot of the geometric information of the space, as they characterize isometries by their boundary action. In this spirit we hope that those cross ratios will be a valuable tool in the studies of symmetric spaces and Euclidean buildings.

We want to refer the reader to section \ref{sec cross ratio preserving maps} to slightly more results in this spirit, e.g. when we get a one-to-one correspondence of Moebius bijections and isometries, and also an analysis in which situations the rescaling of the metric is really necessary.

Concerning the proofs of those theorems: First we show that Moebius bijections split as products of Moebius bijections of irreducible factors; and that Moebius bijections can be extended to building isomorphisms. For rank one symmetric spaces and rank one thick Euclidean buildings it is already known that Moebius bijections extend to isometries. For irreducible thick combinatorial Euclidean buildings it will be enough that Moebius maps are restrictions of building isomorphisms to the chamber set. 
For symmetric spaces and (general) thick Euclidean buildings, we derive additional properties of the building map, using the cross ratio. Those properties will allow us to use theorems (essentially due to Tits) showing that the according maps can be extended to isometries.\medskip

The structure of the paper is as follows. In the preliminaries we recall well known facts of symmetric spaces and Euclidean buildings (we assume the reader to be familiar with those objects) and show basic lemmas we need later on.
In section 3 we define $\R$-valued cross ratios, show basic properties, and illustrate the objects with two examples.
In section 4 we show that the collections of $\R$-valued cross ratios fit together to vector valued cross ratios and motivate that these are the natural objects to consider.
In the last section, section 5, we show that Moebius maps on the chamber set extend to isometries.\medskip 

\textsc{Related Work:} In \cite{Kim} I. Kim constructed a cross ratio very similar to our $\R$-valued cross ratio (Definition \ref{def cross ratio}). Labourie \cite{Labourie-IHES} has given one of the cross ratios in Example \ref{ex Sl n} ad-hoc and used it as tool to understand Hitchin representations. Moreover, Martone and Zhang \cite{Martone-Zhang} have constructed cross ratios on boundaries of surface groups, which in particular for $\text{SL}(n,\R)$-Hitchin representations coincide with the pullback under the boundary map of some of the cross ratios in Example \ref{ex Sl n}.
In \cite{Sambarino-AnnFourier} (see also \cite{Bochi-Potrie-Sambarino}) there is a Gromov product defined, which is closely related to ours.\medskip

\textsc{Acknowledgment:} I want to thank Viktor Schroeder very much for suggesting this topic to me and helping me with fruitful discussions and advice; Linus Kramer for helping me understanding and applying building theory; Beatrice Pozzetti for several helpful comments; and Thibaut Dumont for a valuable comment concerning wall trees. 

\section{Preliminaries}

We use the notation that $M$ is either a symmetric spaces of non-compact type or a thick Euclidean building, $X$ is a symmetric space of non-compact type and $E$ is a thick Euclidean building. In the case of a symmetric space when writing \emph{affine apartment} we mean a maximal flat.

A reference for symmetric spaces of non-compact type is e.g. \cite{Eberlein}; for Euclidean buildings we refer to \cite{Kramer-Weiss}, \cite{Parreau}, \cite{Tits} and also \cite{Kleiner-Leeb}.\footnote{We will use the definition due to \cite{Tits}, which is equivalent to the axioms in \cite{Kramer-Weiss} and \cite{Parreau}, while the definition in \cite{Kleiner-Leeb} would additionally assume metrically completeness.}\medskip

\textsc{Coxeter complex and spherical buildings} \cite{Abramenko-Brown}: Let $W$ be a finite Coxeter group and $S$ the standard set of generators consisting of involutions. Then $W$ can be realized as a reflection group along hyperplanes in $\R^r$ with $r=|S|$. The hyperplanes decompose $\R^r$ and the unit sphere $S^{r-1}$ into (cones over) simplical cells. The maximal, i.e. $r$-dimensional, closed cells in $\R^r$ are called \emph{Weyl sectors}. Lower dimensional cells will be called \emph{conical cells}. The maximal, i.e. $r-1$-dimensional, closed simplical cells in $S^{r-1}$ are called \emph{Weyl chambers}. The set $S$ corresponds to exactly the hyperplanes bounding a Weyl sector. This Weyl sector will be called the \emph{positive sector}, the corresponding chamber in $S^{r-1}$ will be called \emph{positive chamber}.
We can give each simplex adjacent to the positive chamber or positive sector a different label. Then the action of $W$ on the simplical complex induces a unique labeling for all simplices. A fixed label will be called \emph{type}.

In this paper we refer to $(\R^r,W)$ as the \emph{Coxeter complex} and to $(S^{r-1},W)$ as the \emph{spherical Coxeter complex}.

A \emph{spherical building} is a simplical complex $B$ together with a collection of subcomplexes $\mathrm{Apt}(B)$, called \emph{apartments}, which are isomorphic to a fixed spherical Coxeter complex $(S^{r-1},W)$, such that the following holds: 
\begin{enumerate}
\item For any two simplices $a,b\in B$ there is an apartment $A\in \mathrm{Apt}(B)$ with $a,b\in A$
\item If $A,A'$ are apartments containing the simplices $a,b$, then there is a type preserving simplical isomorphism $A\to A'$ fixing $a,b$. 
\end{enumerate}

We say that the building is \emph{modeled over the spherical Coxeter complex $(S^{r-1},W)$}.

A spherical building is called \emph{thick} if each non-maximal simplex is contained in at least three chambers. A (spherical) Coxeter complex is called \emph{irreducible} if the Coxeter group can not be written as a product $W=W_1\times W_2$ of two nontrivial Coxeter groups. A spherical building is called \emph{irreducible} if the spherical Coxeter complex over which it is modeled is irreducible. If a building $B$ is reducible, i.e. modeled over the spherical Coxeter complex $W_1\times W_2$, then it can be written as the spherical join of two buildings, i.e. $B=B_1\circ B_2$ for two spherical buildings $B_1,B_2$ modeled over $W_1,W_2$ respectively and $\circ$ being the spherical join \cite[Sc.3.3]{Kleiner-Leeb}.

Given a simplex $x\in B$ with $B$ a thick spherical building. We denote by $\Res(x):=\{y\in B~|~x\subsetneq y\}$ and call this the \emph{residue of $x$}. Let $A$ be an apartment containing $x$, i.e. a Coxeter complex containing $x$. Let $W$ be the Coxeter group of $A$ and denote by $W_x$ the stabilizer of $x$ under $W$. If $x$ is not a chamber then $\Res(x)$ is itself a spherical building modeled over the Coxeter complex to $W_x$ \cite[3.12]{Tits2}.\medskip

\textsc{Euclidean buildings} \cite{Kramer-Weiss}, \cite{Parreau}, \cite{Tits}, \cite{Kleiner-Leeb}: Let $\hat{W}$ be an affine Coxeter group, i.e. $\hat{W}$ can be realized as a subgroup of the isometry group of $\R^r$ and can be decomposed as a semi-direct product $\hat{W}=W\ltimes T_W$, where $W$ is a finite reflection group and $T_W<\R^r$ is a co-bounded subgroup of translations. Here we assume $r=|S|$, where $S$ is the standard generating set of $W$.
Moreover, let $(E,d)$ be a metric space. A \emph{chart} is an isometric embedding $\phi:\R^r\to E$, and its image is called \emph{affine apartment}; the image of a Weyl sectors and \emph{conical cells} are again called \emph{Weyl sectors} and \emph{conical cells}. Two charts $\phi,\psi$ are called \emph{$\hat{W}$-compatible} if $Y=\phi^{-1}\psi (\R^r)$ is convex in the Euclidean sense and if there is an element $w\in \hat{W}$ such that $\psi\circ w_{|Y}=\phi_{|Y}$. A metric space $E$ together with a collection of charts $\mathcal{C}$, called \emph{apartment system}, is called a \emph{Euclidean building (modeled over the Coxeter group $\hat{W}$)} if it has the following properties:
\begin{enumerate}
\item For all $\phi\in \mathcal{C}$ and $w\in \hat{W}$, the composition $\phi\circ w$ is in $\mathcal{C}$.

\item Any two points $p,q\in E$ are contained in some affine apartment.

\item The charts are $\hat{W}$-compatible.

\item If $a,b\subset E$ are Weyl sectors, then there exists an affine apartment $A$ such that the intersections $A\cap a$ and $A\cap b$ contain Weyl sectors.

\item If $A$ is an affine apartment and $p\in A$ a point, then there is a 1-Lipschitz retraction $\rho:E\to A$ with $d(p,q)=d(p,\rho(q))$ for all $q\in E$.
\end{enumerate}

From this properties it follows that the metric space $E$ is necessarily $\CAT$. The dimension of $\R^r$ is called the \emph{rank} of $E$, i.e. $\mathrm{rk}(E)=r$. While the definition depends on a fixed set of affine apartments, there is always a unique maximal set of affine apartments, called the \emph{complete apartment system}. A set is an affine apartment in the complete apartment system if and only if it is isometric to $\R^r$. In the ongoing \textbf{we will always consider $E$ with its complete apartment system}. If the subgroup of translations $T_W$ is discrete and $E$ is locally compact we call $E$ a \emph{combinatorial} Euclidean building.\medskip

\textsc{Symmetric spaces} \cite[Ch.2]{Eberlein}: Let $X$ be a symmetric space. We will always assume that $X$ is \emph{of non-compact type} and be $d:X\times X\to [0,\infty)$ the natural metric. Moreover, be $G=Iso_0 (X)$, i.e. the connected component of the identity of the isometry group. 

Let $\g= Lie(G)$ and $\g=\kk+\p$ the Cartan decomposition. Fixing a maximal flat $F$ in $X$ together with a basepoint $o\in F$ yields the identification $T_o M\cong\p$. This identification is such that $T_o F\cong \am$ where $\am$ a maximal abelian subspace of $\p$. The restricted root system of $\g$ with respect to $\am$ defines hyperplanes in $\am$ - namely the zero sets of the restricted roots. The Weyl group $W$ of $X$ is the group generated by the reflections along those hyperplanes with respect to the metric that $\am$ inherits from $T_o F\subset T_o X$. Hence we can associate to $X$ a Coxeter complex $(\am,W)$. Let $\am_1$ be the unit sphere in $\am$, then we also get a spherical Coxeter complex $(\am_1,W)$. It is well known that up to isometry the Coxeter complex is independent of the choices. We fix a Weyl sector in $\am$ which we denote by $\am^{+}$ and call \emph{positive sector}. Then $\apo$ will be called the \emph{positive chamber}.\footnote{Usually $\am^{+}$ is called positive Weyl chamber. However, as we will consider Euclidean buildings and symmetric spaces at the same time and we want to distinguish between spherical chambers and cones, we change the usual notation.}
The \emph{rank} of $X$ is the usual rank and equals $\rk(X)=\dim \am$. To keep the notation consistent with buildings we will call maximal flats in $X$ \emph{affine apartments}.\medskip

\textsc{The ideal boundary and Busemann functions} \cite[Part II, Ch.8]{BH}: 
We denote by $\di M$ the ideal boundary; equipped with the cone topology $\di M$ is naturally a topological space.
For every $o\in M$ and every $x\in\di M$ we denote by $\gamma_{ox}$ the unique unit-speed geodesic ray joining $o$ to $x$, i.e. $\gamma_{ox}(0)=o$ and $\gamma_{ox}$ in the class of $x$.
For $o,p,q\in M$ the \emph{Gromov product on $M$} is defined by $(p|q)_o=\frac{1}{2}(d(o,p)+d(o,q)-d(p,q))$.
Let $o\in M$ and $x,y\in \di M$. Then $(\cdot|\cdot)_o:\di M\times \di M\to [0,\infty]$, the \emph{Gromov product} with respect to $o$, is given by $(x|y)_o=\lim\limits_{t\to\infty} (\gamma_{ox}(t)|\gamma_{oy}(t))_o= \lim\limits_{t\to\infty} t-\frac{1}{2} d(\gamma_{ox}(t)|\gamma_{oy}(t))$.
We remark that the convexity of the distance function guarantees the existence of the limit in $[0,\infty]$.

Given $x\in \di M$  the \emph{Busemann function with respect to $x$}, which will be denoted by $b_x:M\times M\to (-\infty,\infty)$, is defined by
\begin{align*}
b_x(o,p)=\lim\limits_{t\to\infty} d(o,\gamma_{px}(t)) - d(p,\gamma_{px}(t))=\lim\limits_{t\to\infty} d(o,\gamma_{px}(t)) - t.
\end{align*}
It holds that $-d(o,p)\leq b_x(o,p)=-b_x(p,o)\leq d(o,p)$ and $b_x(o,p)+b_x(p,q)=b_x(o,q)$ for $o,p,q\in M$. Moreover, it follows directly that $b_x(o,\gamma_{ox}(s))=s$ for all $s\geq 0$ and for all $s\in\R$ if $\gamma_{ox}$ is extended bi-infinitely.

An easy argument in Euclidean geometry yields that the level sets of Busemann functions in $\R^n$ with respect to $x$ in the boundary sphere are affine hyperplanes orthogonal to the direction $x$. In general Busemann level sets with respect to one coordinate are called \emph{horospheres} and the collection of horospheres is independent of the choice of the other coordinate.

The isometry group $Iso(M)$ acts naturally by homeomorphisms on $\di M$, since they map equivalence classes of geodesic rays to equivalence classes of geodesic rays. Moreover, by definition of the Busemann function, it follows $b_x(o,p)=b_{g\cdot x}(g\cdot o,g\cdot p)$ for every $g\in Iso(M)$.\medskip

\textsc{The building at infinity} \cite[Ch.3]{Eberlein}, \cite{Kramer-Weiss}, \cite{Parreau}, \cite{Tits}, \cite{Kleiner-Leeb}: Let $M$ now be either a symmetric space or a Euclidean building. To keep notation simple, we will denote by $(\am,W)$ also the Coxeter complex over which a Euclidean building is modeled. Moreover, $\am_1$ is the unit sphere in $\am$ and hence $(\am_1,W)$ a spherical Coxeter complex. We fix a positive Weyl sector $\am^{+}\subset \am$ and the according positive chamber $\apo=\am_1\cap \am^{+}$. Let $S$ denote the generating set of $W$ consisting of reflections along the walls of $\am^{+}$. By definition we have $\rk(M)=\dim\am$.

The ideal boundary $\di M$ carries naturally the structure of a spherical building $\Di M$ modeled over the spherical Coxeter complex $(\am_1,W)$. The building $\Di M$ will be called the \emph{building at infinity}. 
 
For a Euclidean building $E$ the building at infinity arises as follows: Let $A\subset E$ be an affine apartment. Then $A$ being the image of $(\am,W)$ under a chart implies that $A$ is decomposed into conical cells.
Each conical cell defines a simplex in $\di E$ by taking the geodesic rays contained in the cell for all times. One can show that two conical cells define the same set in $\di E$ if and only if they have finite Hausdorff distance. In the latter case we say the conical cells are equivalent. Taking all conical cells in $E$ modulo the equivalence relation yields a simplical structure on $\di E$; which can be shown the be a spherical building over the spherical Coxeter complex $(\am_1,W)$. 

In a very similar way we get the building at infinity of symmetric spaces $X$: Every maximal flat $F$ with fixed basepoint can be isometrically identified with $\am$. Then the conical cells of $\am$ descend to conical cells in $F\subset X$. Again taking all conical cells in $X$ modulo the equivalence relation of finite Hausdorff distance gives $\di X$ a simplicial structure, which yields a spherical building modeled over $(\am_1,W)$.

Apartments in $\Di M$ correspond to the ideal boundaries of affine apartments of $M$. It is well known that $\Di X$ is a thick building. 
We call an Euclidean building \emph{thick} if in the case $\rk(E)\geq 2$ we have that $\Di E$ is thick, and in the case $\rk(E)=1$ we have that $|\di E|\geq 3$, i.e. $E\neq \R$.

In particular the following important property holds: To every two points $p,q\in M\cup \di M$ we find an affine apartment $A$ in $M$ such that $p,q\in A\cup \di A$. We say that $A$ \emph{joins} $p$ and $q$.

Given two affine apartments $A,A'$ in a Euclidean building $E$ that have a common chamber at infinity, i.e. $c\in \Di E$ such that $c\subset \di A$ and $c\subset\di A'$. Then the intersection $A\cap A'$ contains a Weyl sector with $c$ being its boundary at infinity. Such a Weyl sector is called a \emph{common subsector} of $A$ and $A'$.\medskip

\textsc{The type map} \cite[Sc.4.2.1]{Kleiner-Leeb},\cite[Sc.2.4]{KLP}: To the visual boundary $\di M$ with the building structure $\Di M$ there exists a map $\ty:\di M\to \apo$, called \emph{type map}. Given $x\in\di M$ there is a chamber $c_x\in \Di M$ with $x\in c_x$ and an affine apartment $A$ with $c_x\subset\di A$. Then this yields a isometry from $c_x$ to $\apo$ with respect to the Tits metric on $c_x$ and the angular metric on $\apo$. 
In this way we can assign to each element of $\di M$ an unique element of $\apo$. It can be shown that the image is independent of the chamber and the apartment chosen, hence we get a well defined map $\ty:\di M\to \apo$. The type map is consistent with the types of the spherical building $\Di M$, i.e. two simplices of $\Di M$ are of the same type if and only if they are mapped to the same face of $\apo$ under $\ty$. Hence we also call the faces of $\apo$ \emph{types} ($\apo$ will be a face of itself). 
When speaking of types we denote $\sm=\apo$, i.e. a simplex of $\Di M$ is a chamber if and only if it is of type $\sm$. Faces of $\sm$ will usually be denoted by $\tm$. The set of simplices in $\Di M$ of type $\tm$ will be denoted by $\Ftm(M)$, or just by $\Ftm$ if $M$ is clear out of the context and will be called \emph{flag space}. If we consider chambers we denote this by $\Ftsm$ and call it \emph{full flag space}.

We (ambiguously) call elements in $\xi\in \sm= \apo$ \emph{types}. However, out of the context it is clear if an element or a simplex is meant. We denote by $\inn(\tm)$ the interior of a simplex (and set the interior of a point to be the point itself).
Given a simplex $x\in \Ftm(M)$ and $\xi\in\tm$, we denote by $x_{\xi}$ the unique point in $x\subset\di M$ of type $\xi$.

Let $F:M_1\to M_2$ be an isometry between either two symmetric spaces or two thick Euclidean buildings. Restricting $F$ to the ideal boundary $\di M_1$ induces a building isomorphism $F_{\infty}:\Di M_1\to \Di M_2$. The map $F_{\infty}$ is in general not type preserving. However, that $M_1,M_2$ are isometric implies that they are modeled over the same Coxeter complex and hence have the same fundamental chamber $\sm$. Then we can associate to $F$ a type map $F_{\sm}:\sm\to\sm$ such that $\ty(F_{\infty} (x))=F_{\sm} (\ty(x))$ for every $x\in \di M_1$ and $F_{\sm}$ is a isometry with respect to the angular metric. Moreover, $F(\Ftm(M_1))=\mathrm{Flag}_{F_{\sm}(\tm)}(M_2)$.\medskip

\textsc{The G-action and flag manifolds} \cite[Ch.3]{Eberlein}, \cite[Sc.2.4]{KLP}: Let $X$ be a symmetric space and $G=Iso_0(X)$. Then the cone topology on $\di X$ induces a topology on $\Di X$ such that all flag spaces are compact. Moreover, given $x\in\Ftm(X)$, let $P_x$ denote the stabilizer of $x$ under the $G$-action. Then we can identify $\Ftm(X)\simeq G\slash P_x$ with the identification being $G$-equivariant and homeomorphic; the group $P_x$ is a parabolic subgroup of $G$ and $G\slash P_x$ is equipped with the quotient topology of the topological group $G$. 
Moreover, $\Ftm(X)\simeq G\slash P_x$ yields a smooth structure on $\Ftm(X)$ (inherited from $G\slash P_x$) making it a compact connected manifold. The spaces $G\slash P_x$ are called \emph{Furstenberg boundaries} or \emph{flag manifolds} (motivating our notion of flag space).
Let $K$ be a maximal compact subgroup of $G$. Then already $K$ acts transitive on the flag manifolds and given $x\in\Ftm(X)$ we can identify $\Ftm(X)\simeq K\slash K_x$ $K$-equivariant and homeomorphically, where $K_x=\st_K (x)$.
Moreover, we remark that the $G$-action is type preserving, i.e. $g_{\sm}=\id$ for all $g\in G$.\medskip

\textsc{The opposition involution:} An important map for us will be the \emph{opposition involution} $\iota:\am\to \am$, which is given by $\iota=-\id\circ w_0$ with $w_0\in W$ the maximal element of the Coxeter group with respect to the generating set $S$. If $W$ is an irreducible Weyl group, then $\iota=id$ if and only if $W$ is not of type $A_n$ with $n\geq 2$, $D_{2n+1}$ with $n\geq 2$ or $E_6$ \cite[2.39]{Tits2}. Moreover, we remark that we can restrict $\iota:\am^{+}_1\to\am^{+}_1$ and that $\iota$ is an isometry with respect to the angular metric.\medskip

\textsc{Opposite simplices} \cite[Sc.2.2,\:2.4]{KLP}: There is a natural notion of \emph{opposition} in spherical buildings. This corresponds to the following: Let $x,y\in \Di M$ and let $A_{\infty}$ be an apartment in $\Di M$ such that $x,y\in A_{\infty}$. Since $A_{\infty}$ can be identified with the unit sphere $\am_1$, there is a natural map $-\id: A_{\infty}\to A_{\infty}$. 
Then $x$ is \emph{opposite} of $y$, denoted by $x\op y$, if and only if $x=-\id (y)$. The action of the spherical Coxeter group $W$ leaves the type invariant. Therefore, assume for the moment that $W$ is modeled in $A_{\infty}$ and $x$ is a face of the positive chamber. Denote by $w_0: A_{\infty}\to A_{\infty}$ the maximal element of $W$. Then $w_0(y)$ is a face of the positive chamber and of the same type as $y$ and hence $y$ is of type $-\id\circ w_0 (x)=\iota x$. Hence all simplices opposite of elements in $\Ftm$ are contained in $\Ftim$.
For later use we denote
\begin{align}\label{eq def domain of definition} 
  &\Atop:=\{ (x_1,y_1,x_2,y_2)\in (\Ftm\times \Ftim)^2~|~x_1,x_2\op y_1,y_2 \}\\
 \At:=\{ (x_1 &,y_1,x_2,y_2)\in (\Ftm\times \Ftim)^2~|~x_i\op y_i \text{   or   }x_i\op y_j , i,j=1,2, i\neq j \}. \nonumber
\end{align}

Opposition of simplices has the following important connection to bi-infinite geodesics: Let $z_1,z_2\in \di M$ and $A\subset M$ an affine apartment with $z_1,z_2\in\di A$. Then one can show that there exists a bi-infinite geodesics joining $z_1$ and $z_2$ if and only if there exists one in $A$. From Euclidean geometry it follows that the $z_i$ can be joined by a bi-infinite geodesic in $A$ if and only if $z_1=-\id(z_2)$ with $-\id: \di A\to \di A$ as before. This can easily be seen to be equivalent to the unique simplices $\tm_{z_i}\in \Di M$ containing the $z_i$ in its interior being opposite, i.e. $\tm_{z_1}\op \tm_{z_2}$, and $\ty(z_1)=\iota \ty(z_2)$. 

We will call points $z_1,z_2\in\di M$ \emph{opposite} if they can be joined by a bi-infinite geodesic and denote this also by $z_1\op z_2$. 
Moreover,  for every $\xi\in \tm$ and $(x,y)\in \Ftm\times \Ftim$ with $x\op y$, it follows that $x_{\xi}$ is opposite to $y_{\iota\xi}$.\medskip

\textsc{Symmetric spaces, Langlands decomposition \cite[Sc.2.17]{Eberlein}, \cite[Sc.2.10]{KLP}:} In case of a symmetric space $X$ and given $x\in\Ftm(X)$, the set of simplices opposite to $x$ is an open and dense subset of $\Ftim(X)$ (which can be deduced from the Bruhat decomposition of $G\slash P$). Moreover, for $(x,y)\in \Ftm(X)\times \Ftim(X)$ we have $x\op y$ if and only if the pair is in the unique open and dense $G$-orbit in $\Ftm(X)\times \Ftim(X)$.
In particular, it follows in this case that $\At$ and $\Atop$ are open and dense subsets of $(\Ftm\times \Ftim)^2$.

Every parabolic subgroup $P_x$ has a natural decomposition $P_x=K_x A_x N_x$ called the \emph{Langlands decomposition}, where $K_x$ is compact and $N_x$ is nilpotent. The group $N_x$ is called \emph{horospherical subgroup} and is unique, while $K_x$ and $A_x$ are not. 
The horospherical subgroup has several important properties; it leaves the Busemann function with respect to $x_{\xi}\in x\in\Ftm(X)$ invariant, i.e. $b_{x_{\xi}}(o,p)=b_{x_{\xi}}(n\cdot o,p)=b_{x_{\xi}}(o,n\cdot p)$ for all $n\in N_x$ and $\xi\in\tm$; given a geodesic ray $\gamma_{x_{\xi}}$ with endpoint in $x\subset\di X$, we have $d(\gamma_{x_{\xi}}(t),n\cdot \gamma_{x_{\xi}}(t))\to 0$ for $t\to\infty$ for all $n\in N_x$;
moreover, $N_x$ acts simply transitive on the set of simplices opposite to $x$. If $x$ is a chamber, i.e. $x\in \Ftsm(M)$, then $N_x$ acts simply transitive on the set of maximal flats containing $x$ in its boundary.\medskip

\textsc{Parallel sets \cite[Sc.2.11, 2.20]{Eberlein},\cite[Sc.2.4]{KLP},\cite[Sc.4.8]{Kleiner-Leeb}:} Let $(x,y)\in \Ftm(M)\times \Ftim(M)$ with $x\op y$ and let $\xi$ be an element of $\inn(\tm)$. Then the \emph{parallel set} with respect to $x,y$, denoted by $P(x,y)$, is the set of all points that lie on a bi-infinite geodesic joining $x_{\xi}$ to $y_{\iota\xi}$.

The parallel sets split metrically as products, i.e. $P(x,y)\simeq F_{xy}\times CS(x,y)$, where $F_{xy}$ is an isometrically embedded $\R^n$ such that $x,y\subset \di F_{xy}$ and $x,y$ are simplices of maximal dimension in the sphere $\di F_{xy}$ - in particular the dimension of the spherical simplices $x,y$ equals $n-1$. Then it follows that the parallel set is independent of the choice of type $\xi\in\inn(\tm)$, as for each type $\xi\in\inn(\tm)$ geodesics in $M$ joining $x_{\xi},y_{\iota\xi}$ are of the form $(\gamma_{x_{\xi}y_{\iota\xi}}(t),p)$ with $\gamma_{x_{\xi}y_{\iota\xi}}$ a geodesic in $F_{xy}$ joining $x_{\xi},y_{\iota\xi}$ and $p$ is a point in $CS(x,y)$.
 
The space $CS(x,y)$ is called \emph{cross section}. In case of a symmetric space $X$ the cross section is itself a symmetric space without Euclidean de Rham factors, in case of a Euclidean building the cross section is again a Euclidean building. In both cases the rank is given by $\mathrm{rk}(CS(x,y))=\mathrm{rk}(M)-\dim F_{xy}$

Let $\tm$ be a face of $\sm=\am_1$. Let $\at$ be the subspace of $\am$ defined by $\tm$, i.e. the smallest subspace of $\am$ containing $\tm$ and $0$. Let $\xi_1,\ldots,\xi_k\in\am$ be the corners of the spherical simplex $\tm$. Then $\at=\mathrm{span}_{i=1,\ldots,k} ~ \xi_i$. It is immediate that we can also identify $P(x,y)\simeq \at\times CS(x,y)$. We can additionally impose that this identification is in such a way that $x\simeq \di \at^{+}$ where $\at^{+}:=(\at\cap \ap)$.

\bl\label{lem horospheres in parallel set}
Let $(x,y)\in\Ftm\times\Ftim$ with $x\op y$ and be $p,q\in P(x,y)$. Let $\pi: P(x,y)\simeq \at\times CS(x,y)\to\at$ be the projection to the first factor. Then for each $\xi\in\tm$ we have that $b_{x_{\xi}}(p,q)=(b_{x_{\xi}})_{|\at}(\pi(p),\pi(q))$, i.e. the Busemann function is independent of the second factor of the product.
\el

\begin{proof}
Let $\gamma_{qx_{\xi}}$ denote the geodesic ray from $q$ to $x_{\xi}$. Moreover, be $q=(q_1,q_2)$ under the identification $P(x,y)\simeq \at\times CS(x,y)$. Then we have that $\gamma_{qx_{\xi}}\simeq (\gamma_{q_1 x_{\xi}},q_2)$ where $\gamma_{q_1 x_{\xi}}$ is the geodesic ray in $\at$ from $q_1$ to $x_{\xi}$. 
Using that metrically $P(x,y)\simeq \at\times CS(x,y)$ and $p=(p_1,p_2)$ we derive $d(p,\gamma_{qx_{\xi}}(t))=\sqrt{d(p_1,\gamma_{q_1x_{\xi}}(t))^2+d(p_2,q_2)^2)}$. If we set $K_2:=d(p_2,q_2)^2$, then $b_{x_{\xi}}(p,q)=\lim_{t\to\infty} \sqrt{d(p_1,\gamma_{q_1x_{\xi}}(t))^2+K_2}-t$. As $p_1, \gamma_{q_1x_{\xi}}(t)\in\at$, it reduces to Euclidean geometry, i.e. $d(p_1,\gamma_{q_1x_{\xi}}(t))=\sqrt{b_{x_{\xi}}(p_1,\gamma_{q_1x_{\xi}}(t))^2+K_1}$ with $K_1$ the squared distance from $p_1$ to the (now) bi-infinite geodesic $\gamma_{q_1x_{\xi}}$. It follows that $b_{x_{\xi}}(p_1,\gamma_{q_1x_{\xi}}(t))=t+b_{x_{\xi}}(p_1,q_1)$. Using a substitution $t=s^{-1}$ and a Taylor series for the root expression below yields
\begin{align*}
 b_{x_{\xi}}(p,q)=&\lim_{t\to\infty} \sqrt{(t+b_{x_{\xi}}(p_1,q_1))^2+K_1+K_2}-t\\
=&\lim_{s\to 0} s^{-1}(\sqrt{(1+2sb_{x_{\xi}}(p_1,q_1)+s^2(b_{x_{\xi}}(p_1,q_1)^2 +K_1+K_2)}-1)\\
=&
b_{x_{\xi}}(p_1,q_1).\qedhere
\end{align*}
\end{proof}

We will also need the following lemma.

\bl\label{lem busemann in parallel set} Let $(x,y)\in \Ftm\times \Ftim$ with $x\op y$ and $\xi\in \tm$. Moreover let $p_1,p_2\in P(x,y)$. Then $b_{x_{\xi}}(p_1,p_2)=-b_{y_{\iota\xi}}(p_1,p_2)$.
\el

\begin{proof}
Let $\gamma_i,\; i=1,2$ be bi-infinite geodesics with $\gamma_i(0)=p_i$, $\gamma_i(+\infty)=x_{\xi}$ and $\gamma_i(-\infty)=y_{\iota\xi}$, which exists by assumption. The $\gamma_i$ are parallel and denote by $C$ their distance. Then the Flat Strip Theorem (see e.g. \cite{BH}) implies that the convex hull of $\gamma_1(\R)\cup \gamma_2(\R)$ is isometric to a flat strip $\R \times [0,C] \subset \R^2$ with $\gamma_i$ identified with $\R\times 0$, $\R\times C$ respectively.

It follows that the level sets of the Busemann function $b_{x_{\xi}}(\cdot,p_2)$ in $\R\times [0,C]$ are given by hyperplanes orthogonal to $\gamma_i$, i.e. are of the form $s\times [0,C]$ and the same holds for $b_{y_{\iota\xi}}(\cdot,p_2)$. In addition, $\gamma_i$ joining $x_{\xi}$ to $y_{\iota\xi}$ implies $b_{x_{\xi}}(\cdot,p_2)_{|\gamma_i}=-b_{y_{\iota\xi}}(\cdot,p_2)_{|\gamma_i}$. Then the claim is direct consequence.
\end{proof}

\textsc{Retracts \cite{Parreau}:} Lastly, we need to introduce the notion of retracts of $M$ to affine apartments with respect to chambers at infinity. For the construction we will distinguish between Euclidean buildings and symmetric spaces.

Let $E$ be a Euclidean building. Let $A\subset E$ be an affine apartment and $x\subset \di A$ a chamber of the building at infinity. Then there exists a 1-Lipschitz map $\rho_{x,A}:E\to A$ which is an isometry when restricted to any affine apartment $A'$ with $x\subset \di A'$ (i.e. any affine apartment that contains the chamber $x$ in its boundary), and the identity on $A$ \cite[Prop.1.20]{Parreau}. We call this map \emph{(horospherical) retract} with respect to $x$. Horospherical retracts have the following important property:

\bl Let $\rho_{x,A}:E\to A$ be a horospherical retract with respect to $x\in \Ftsm(E)$. Then $b_{x_{\xi}}(o,p)=b_{x_{\xi}}(\rho_{x,A}(o),p)=b_{x_{\xi}}(o,\rho_{x,A}(p))$ for all $o,p\in E$ and $\xi\in\sm$.
\el

\begin{proof}
To $o\in E$ there exists an affine apartment $A_o$ containing $o$ and $x\subset \di A_o$. As mentioned, the horopsheres with respect to $x_{\xi}$ in $A_o$ are hyperplanes orthogonal to the direction $x_{\xi}$.

By construction, the two affine apartments $A$, $A_o$ have the same chamber in its boundary, which implies that they have a common subsector. Hence $\rho_{x,A}$ is the identity on the non-empty intersection $A\cap A_o$. Moreover, $\rho_{x,A}$ is an isometry when restricted to $A_o$. Since $\rho_{x,A}$ leaves each horosphere intersecting $A\cap A_o$ invariant, it has to map the level set of $b_{x_{\xi}}(\cdot,p)$ in $A_o$ to the corresponding level set in $A$.
The other equality follows for example form the symmetry $b_{x_{\xi}}(o,p)=-b_{x_{\xi}}(p,o)$
\end{proof}

Let $X$ be a symmetric space, $A\subset X$ be a maximal flat (an affine apartment for us) and $x\subset \di A$ a chamber at infinity. To any $o\in X$ there exists a unique maximal flat $A_o$ with $o\in A_o$ and $x\subset \di A_o$. Then we define $\rho_{x,A}(o):=n_{x,A_o}\cdot o$ for $n_{x,A_o}$ the unique element in $N_x$ that maps $A_o$ to $A$. Again we call $\rho_{x,A}:X\to A$ \emph{(horospherical) retract}.\medskip

For later reference: To every affine apartment $A\subset M$ and a chamber $x\subset \di A$ we have a well defined map $\rho_{x,A}:M\to A$ such that
\begin{align}\label{eq horospherical retract and busemann}
b_{x_{\xi}}(o,p)=b_{x_{\xi}}(\rho_{x,A}(o),p)=b_{x_{\xi}}(o,\rho_{x,A}(p))
\end{align}
for all $o,p\in M$ and $\xi\in\sm$. Moreover, it is known that two opposite chambers $x,y\in \Ftsm$ are contained in an unique apartment $A_{\infty}$ of $\Di M$ and this corresponds to an unique affine apartment $A_{xy}\subset M$. Hence to $x,y\in \Ftsm$ with $x\op y$ we set $\rho_{x,y}:=\rho_{x,A_{xy}}$.

\bl \label{lem retracts of geodesics are geos again}
Let $x,y\in\Ftm$ with $x\op y$ and $o\in M$. Then for all $\xi \in\tm$ we have that $\rho_{c_x,c_y}(\gamma_{ox_{\xi}}(t))$ is a geodesic in $P(x,y)$, where $c_x,c_y\in \Ftsm$ such that $x$ is a face of $c_x$, $y$ is a face of $c_y$ and $c_x\op c_y$.
\el

We remark that $x\op y$ implies that such $c_x,c_y\in \Ftsm$ always exist. Namely, take an apartment containing $x$ and $y$. Take $c_x\in\Ftsm$ such that $x$ is a face of $c_x$. Take $c_y\in\Ftsm$ the unique opposite chamber in the apartment. Then $x\op y$ implies that $y$ is a face of $c_y$. 

\begin{proof}
For a symmetric space $X$ this follows since $\rho_{c_x,c_y}$ is the same element of $G$ for all points $\gamma_{ox_{\xi}}(t)$ and that $G<Iso(X)$. Hence $\rho_{c_x,c_y}(\gamma_{ox_{\xi}}(t))$ is the image of a geodesic under an isometry. The image $\rho_{c_x,c_y}(\gamma_{ox_{\xi}}(t))$ is geodesic ray with endpoint $x_{\xi}$ in an affine apartment joining $x$ and $y$. Then $y\op x$ implies that if we extend $\rho_{c_x,c_y}(\gamma_{ox_{\xi}}(t))$ bi-infinitely it joins $x_{\xi}$ to $y_{\iota\xi}$, i.e. this geodesic is contained in $P(x,y)$.

Consider a Euclidean building $E$. Denote by $A_{xy}$ the unique affine apartment joining $c_x$ and $c_y$.
Let $A$ be an affine apartment containing $o$ and $c_x\subset \di A$. Then it follows that $\gamma_{ox_{\xi}}(t)\in A$ for all $t\in \R_{+}$. As $\rho_{c_x,c_y}$ is an isometry on affine apartments containing $c_x$, it follows that $\rho_{c_x,c_y}(\gamma_{ox_{\xi}}(t))\subset A_{xy}$ is the image of a geodesic under an isometry. Since one of the endpoints is $x_{\xi}$, we can extend the geodesic in $A_{xy}$ uniquely to a bi-infinite geodesic joining $x_{\xi}$ and $y_{\iota\xi}$. Thus $\rho_{c_x,c_y}(\gamma_{ox_{\xi}}(t))\subset P(x,y)$.
\end{proof}

\section{Cross ratios}

Let $M$ be a symmetric space of non-compact type or a thick Euclidean building. Let $\sm$ be the fundamental chamber of the associated spherical Coxeter complex and $\tm$ a face of $\sm$. For any type $\xi\in \sm$ such that $\xi\in\inn(\tm)$ and any $o\in M$ we define a \emph{Gromov product $(\,\cdot\,|\,\cdot\,)_{o,\xi}: \Ftm(M)\times \Ftim(M)\to [0,\infty]$ with base-point $o$} by 
\begin{align*}
(x|y)_{o,\xi} := \lim_{t\to\infty} t - \frac{1}{2} d(\gamma_{ox_{\xi}}(t),\gamma_{oy_{\iota\xi}}(t))
\end{align*}
for $(x,y)\in \Ftm(M)\times \Ftim(M)$ and $\gamma_{ox_{\xi}}(t),\gamma_{oy_{\iota\xi}}(t)$ the unit speed geodesics from $o$ to $x_{\xi},y_{\iota\xi}$, respectively. 
Using this we define the \emph{(additive) cross ratio $\cro_{o,\xi}:\At\to [-\infty,\infty]$ with respect to $(o,\xi)$} by
\begin{align*} 
\cro_{o,\xi}(x_1,y_1,x_2, y_2):= -(x_1|y_1)_{o,\xi}-(x_2|y_2)_{o,\xi}+(x_1|y_2)_{o,\xi}+(x_2|y_1)_{o,\xi}
\end{align*}
where $\At$ is the set of quadrupels $(x_1,y_1,x_2,y_2)\subset (\Ftm(M)\times \Ftim(M))^2$ as in equation \eqref{eq def domain of definition}. If $\xi\in\inn(\tm)$, we also denote $\Af:=\At$.
By definition $\cro_{o,\xi}$ has the following symmetries, whenever all factors are defined,
\begin{align}\label{eq symmetries of cross ratios}
\cro_{o,\xi}(x_1,y_1,x_2, y_2) &=-\cro_{o,\xi}(x_1,y_2,x_2, y_1)=-\cro_{o,\xi}(x_2,y_1,x_1, y_2) \nonumber\\
\cro_{o,\xi}(x_1,y_1,x_2, y_2) &=\cro_{o,\xi}(x_1,y_1,w, y_2)\: +\cro_{o,\xi}(w,y_1,x_2, y_2) \\
\cro_{o,\xi}(x_1,y_1,x_2, y_2) &=\cro_ {o,\xi}(x_1,y_1,x_2, v)\: +\cro_{o,\xi}(x_1,v,x_2, y_2). \nonumber
\end{align}
The last two symmetries are called \emph{cocycle identities}.\medskip

\textbf{Notation:} Let $\tm$ be face of $\sm$ and be $\xi\in\partial\tm$. Then we drop for any $(x,y)\in\Ftm\times\Ftim$ the projection maps in the Gromov product (and in the cross ratio) for notational reasons, i.e. $(x|y)_{o,\xi}:=(\pi_{\xi}(x),\pi_{\iota\xi}(y))_{o,\xi}$, where $\tm_{\xi}$ is the face of $\tm$ containing $\xi$ in its interior and $\pi_{\xi}:\Ftm\to \mathrm{Flag}_{\tm_{\xi}}$, $\pi_{\iota\xi}:\Ftim\to \mathrm{Flag}_{\iota\tm_{\xi}}$ are the obvious projection maps. 

\bp\label{prop Gromov product as Busemann function} Let $M$ be a symmetric space or thick Euclidean building, $o\in M$, $(x,y)\in \Ftm(M)\times \Ftim(M)$ with $x\op y$ and $c_x,c_y\in \Ftsm(M)$ such that $x$ is a face of $c_x$, $y$ is a face of $c_y$ and $c_x\op c_y$. Then for every $\xi\in\tm$
\begin{align*}
(x|y)_{o,\xi} =\frac{1}{2}b_{x_{\xi}}(o,\rho_{c_y,c_x}(o))=\frac{1}{2}b_{y_{\iota\xi}}(o,\rho_{c_x,c_y}(o)).
\end{align*}
\ep

\begin{proof}
In case of a symmetric space let $N_x$ be the horospherical subgroup of $P_x=\st(x)$ and be $n_x(o,y)\in N_x$ the unique element such that $n_x(o,y)\cdot o\in P(x,y)$: Extend $\gamma_{ox}$ bi-infinitely and let $z\in\Ftim$ be such that $\gamma_{ox}(-\infty)\in z$. Then $n_x(o,y)\in N_x$ is the unique element with $n_x(o,y)(z)=y$. By construction we have $n_x(o,y)\cdot o\in P(x,y)$. 

We define in the same way $n_y(o,x)\in N_y$ and set $\gamma_{xy}(t):=n_x(o,y)\cdot \gamma_{ox_{\xi}}(t)$ and $\gamma_{yx}(t):=n_x(o,y)\cdot\gamma_{oy_{\iota\xi}}(t)$. Then $\gamma_{xy}, \gamma_{yx}$ are geodesics in $P(x,y)$ with the same (un-ordered) end points. Hence they are parallel. Moreover, $n_x(o,y)\in N_x$ implies that $d(\gamma_{ox_{\xi}}(t),\gamma_{xy}(t))\to 0$ for $t\to \infty$ and similarly $d(\gamma_{oy_{\iota\xi}}(t),\gamma_{yx}(t))\to 0$.

The triangle inequality yields that $(x|y)_{o,\xi} =\lim_{t\to\infty} t-\frac{1}{2} d(\gamma_{xy}(t),\gamma_{yx}(t))$. By construction $\gamma_{xy},\gamma_{yx}$ are parallel geodesics; hence by the Flat Strip Theorem (see e.g. \cite{BH}) the distance $d(\gamma_{xy}(t),\gamma_{yx}(t))$ decomposes into a part parallel to the geodesics and the distance of the images of the geodesics, which is a constant and will be denoted by $C$. 

The part parallel to the geodesics is $b_{x_{\xi}}(\gamma_{yx}(t),\gamma_{xy}(t))$ - or in the same way $b_{y_{\iota\xi}}(\gamma_{xy}(t),\gamma_{yx}(t))$. Using that we have geodesics asymptotic to $x_{\xi}$ we derive that $b_{x_{\xi}}(\gamma_{yx}(t),\gamma_{xy}(t)))=2t + b_{x_{\xi}}(\gamma_{yx}(0),\gamma_{xy}(0))$. Altogether \begin{align}\label{eq proof of Prop 1}
(x|y)_{o,\xi}= &\lim_{t\to\infty} t-\frac{1}{2}d(\gamma_{ox_{\xi}}(t),\gamma_{oy_{\iota\xi}}(t))=\lim_{t\to\infty} t-\frac{1}{2}d(\gamma_{xy}(t),\gamma_{yx}(t)) \nonumber\\
= & \lim_{t\to\infty} t- \frac{1}{2}(\sqrt{(2t+b_{x_{\xi}}(\gamma_{yx}(0),\gamma_{xy}(0)))^2+C^2})\\
= & - \frac{1}{2}b_{x_{\xi}}(\gamma_{yx}(0),\gamma_{xy}(0))=\frac{1}{2}b_{x_{\xi}}(\gamma_{xy}(0),\gamma_{yx}(0)),\nonumber
\end{align}
while the second to last equality follows using Taylor series at $s=0$ after substituting $s=t^{-1}$ (see also the calculations in example \ref{exam h2 cross h2}).\medskip

In case of a Euclidean building $E$, let $A_o$ be an affine apartment containing $\gamma_{ox_{\xi}}(t)$, let $d_x\in\Ftsm$ be such that $d_x\subset \di A_o$ and $x\subset d_x$. Moreover, be $d_y\in\Ftsm$ a chamber opposite to $d_x$ such that $y$ is a face of $d_y$ and let $A_{xy}$ be the unique affine apartment that $d_x$ and $d_y$ define. 

Then the affine apartments $A_o$ and $A_{xy}$ have a common subsector. Hence there exists $T_x\geq 0$ such that for $t\geq T_x$ the geodesic $\gamma_{ox_{\xi}}(t)$ is parallel to a geodesic $\gamma_{xy}$ in the subsector - denote the distance of the geodesic rays by $C_x$; Extend $\gamma_{xy}$ bi-infinite in $A_{xy}$ such that it is in the same horosphere with respect to $x_{\xi}$ as $\gamma_{ox_{\xi}}(t)$ for all (positive) time. That $\gamma_{xy}$ is in $A_{xy}$ with one endpoint being $x_{\xi}$ implies that $\gamma_{xy}$ joins $x_{\xi}$ and $y_{\iota\xi}$ and hence $\gamma_{xy}\subset P(x,y)$.

In the same way we construct $\gamma_{yx}\subset P(x,y)$ to $\gamma_{oy_{\iota\xi}}$ such that those geodesics are parallel for $t\geq T_y$ - denote the distance by $C_y$. Since $\gamma_{xy},\gamma_{yx}$ join the same points at infinity, they are parallel - denote the distance by $C_0$. Then the triangle inequality together with the Flat Strip theorem yields for $t\geq \max\{T_x,T_y\}$ that $d(\gamma_{ox_{\xi}}(2t),\gamma_{oy_{\iota\xi}}(2t))$ is smaller or equal than
\begin{align*}
&d(\gamma_{ox_{\xi}}(2t),\gamma_{xy}(t))+d(\gamma_{xy}(t),\gamma_{yx}(t)) +d(\gamma_{yx}(t),\gamma_{oy_{\iota\xi}}(2t))\\
=& \sqrt{t^2-C^2_x}+\sqrt{b_{x_{\xi}}(\gamma_{yx}(t),\gamma_{xy}(t))^2+C_0^2}+\sqrt{t^2-C^2_y}
\end{align*}
Since $\gamma_{xy}$ and $\gamma_{yx}$ are asymptotic to $x_{\xi}$, we derive that $b_{x_{\xi}}(\gamma_{yx}(t),\gamma_{xy}(t)))=2t + b_{x_{\xi}}(\gamma_{yx}(0),\gamma_{xy}(0))$. Therefore
\begin{align*}
(x|y)_{o,\xi}\geq \lim_{t\to\infty} 2t-\frac{1}{2}(\sqrt{t^2-C_x^2}+\sqrt{(2t + b_{x_{\xi}}(\gamma_{yx}(0),\gamma_{xy}(0)))^2+C_0^2}+\sqrt{t^2-C_y^2}).
\end{align*}
We substitute $t=s^{-1}$. Then a Taylor expansions for the root expressions at $s=0$ yields that $(x|y)_{o,\xi}\geq -\frac{1}{2} b_{x_{\xi}}(\gamma_{yx}(0),\gamma_{xy}(0))=\frac{1}{2} b_{x_{\xi}}(\gamma_{xy}(0),\gamma_{yx}(0))$. 

We claim that $\lim_{t\to\infty} b_{x_{\xi}}(\gamma_{yx}(t),\gamma_{xy}(t)))-b_{x_{\xi}}(\gamma_{oy_{\iota\xi}}(t),\gamma_{ox_{\xi}}(t))=0$: 
By construction $b_{x_{\xi}}(\gamma_{xy}(t),\gamma_{ox_{\xi}}(t))=0$. Therefore it is enough to show that $\lim_{t\to\infty}b_{x_{\xi}}(\gamma_{oy_{\iota\xi}}(t),\gamma_{yx}(t))=0$, as Busemann functions satisfy $b_z(p,q)+b_z(q,o)=b_z(p,o)$.

By construction we have that the geodesic $\gamma_{yx}$ joins $x_{\xi}$ and $y_{\iota\xi}$. Therefore $b_{x_{\xi}}(\gamma_{oy_{\iota\xi}}(t),\gamma_{yx}(t))=\lim_{s\to\infty} d(\gamma_{oy_{\iota\xi}}(t),\gamma_{yx}(t-s))- s$. Moreover,
\begin{align*}
d(\gamma_{oy_{\iota\xi}}(t),\gamma_{yx}(t-s))\leq d(\gamma_{oy_{\iota\xi}}(t),\gamma_{yx}(T_y)) + |t-s-T_y|
\end{align*}
Applying the Flat Strip Theorem with an according Taylor expansion as before, we derive that $\lim_{t\to\infty} d(\gamma_{oy_{\iota\xi}}(t),\gamma_{yx}(T_y))-t\to -T_y$. In particular
\begin{align*}
\lim_{t\to\infty}b_{x_{\xi}}(\gamma_{oy_{\iota\xi}}(t),\gamma_{yx}(t))\leq & \lim_{t\to\infty}(\lim_{s\to\infty} d(\gamma_{oy_{\iota\xi}}(t),\gamma_{yx}(T_y))-t+ s+T_y-s)=0.
\end{align*}

It follows from the definition of Busemann functions that if $q\in M$ lies on a bi-infinite geodesics joining $z,w\in \di M$, then $b_z(p,q)+b_w(p,q)\geq 0$. Hence we derive $b_{x_{\xi}}(\gamma_{oy_{\iota\xi}}(t),\gamma_{yx}(t))+b_{y_{\iota\xi}}(\gamma_{oy_{\iota\xi}}(t),\gamma_{yx}(t))\geq 0$. By construction $b_{y_{\iota\xi}}(\gamma_{yx}(t),\gamma_{oy_{\iota\xi}}(t))=0$. Thus $b_{x_{\xi}}(\gamma_{oy_{\iota\xi}}(t),\gamma_{yx}(t))\geq 0$, which yields the claim.

We have $d(\gamma_{oy_{\iota\xi}}(t),\gamma_{ox_{\xi}}(t)) \geq b_{x_{\xi}}(\gamma_{oy_{\iota\xi}}(t),\gamma_{ox_{\xi}}(t))\to b_{x_{\xi}}(\gamma_{yx}(t),\gamma_{xy}(t))$, for $t\to\infty$. Thus
\begin{align*}
(x|y)_{o,\xi}\leq \lim_{t\to \infty} t - \frac{1}{2}b_{x_{\xi}}(\gamma_{yx}(t),\gamma_{xy}(t))=\frac{1}{2} b_{x_{\xi}}(\gamma_{xy}(0),\gamma_{yx}(0))
\end{align*}
Altogether $(x|y)_{o,\xi}=\frac{1}{2} b_{x_{\xi}}(\gamma_{xy}(0),\gamma_{yx}(0))$.\medskip

Consider a symmetric space or a Euclidean building $M$ and let $\gamma_{xy},\gamma_{yx}$ be the accordingly constructed geodesics. Then $b_{x_{\xi}}(\gamma_{xy}(0),\gamma_{ox_{\xi}}(0))=0$ while $\gamma_{ox_{\xi}}(0)=o$ and also $b_{y_{\iota\xi}}(\gamma_{yx}(0),o)=0$. For notational reasons set $\rho_x:=\rho_{c_x,c_y}$ and $\rho_y:=\rho_{c_y,c_x}$ Then $\rho_y(o),\gamma_{yx}(0)\in P(x,y)$. Together with equation \eqref{eq horospherical retract and busemann} and Lemma \ref{lem busemann in parallel set} this yields
\begin{align*}
b_{x_{\xi}}(\gamma_{xy}(0),\gamma_{yx}(0))= & b_{x_{\xi}}(\gamma_{xy}(0),\rho_x(o))+b_{x_{\xi}}(\rho_x(o),\rho_y(o))+b_{x_{\xi}}(\rho_y(o),\gamma_{yx}(0))\\
= & b_{x_{\xi}}(o,\rho_y(o))-b_{y_{\iota\xi}}(\rho_y(o),\gamma_{yx}(0))=b_{x_{\xi}}(o,\rho_y(o)).
\end{align*}
In a similar way it follows also $b_{x_{\xi}}(\gamma_{xy}(0),\gamma_{yx}(0))=b_{y_{\iota\xi}}(o,\rho_x(o))$. Finally, $(x|y)_{o,\xi}=\frac{1}{2} b_{x_{\xi}}(\gamma_{xy}(0),\gamma_{yx}(0))$ implies the claim.
\end{proof}

\bc \label{corollary-not-opposite-equivalent-distance-zero}
Let $(x,y)\in \Ftm\times\Ftim$ and $o\in M$. Then $(x|y)_{o,\xi}= \infty\Longleftrightarrow x\nop y$. 
\ec

\begin{proof}
Let $(x,y)\in \Ftm\times\Ftim$ be such that $x\nop y$. Let $A$ be an affine apartment containing $x,y$ in its boundary. Let $p\in A$ and $\gamma_{px_{\xi}}$, $\gamma_{py_{\iota\xi}}$ be the unit speed geodesics joining $p$ to $x_{\xi},y_{\iota\xi}$, respectively. A straight forward argument in Euclidean geometry yields that $d(\gamma_{px_{\xi}}(t), \gamma_{py_{\iota\xi}}(t))= 2 \alpha t$ with $\alpha$ depending on the angle of the geodesics. Then $x\nop y$ implies that $\gamma_{px_{\xi}}(t)\neq \gamma_{py_{\iota\xi}}(-t)$ and hence $\alpha<1$, i.e. $(x|y)_{p,\xi}= \infty$.

Now let $\gamma_{ox_{\xi}}, \gamma_{oy_{\iota\xi}}$ be the unit speed geodesics joining $o$ to $x_{\xi},y_{\iota\xi}$, respectively. Since $\gamma_{ox_{\xi}}$ and $\gamma_{px_{\xi}}$ define the same point in the ideal boundary, we can derive - by the convexity of the distance functions along geodesics in non-positive curvature - that $d(\gamma_{ox_{\xi}}(t),\gamma_{px_{\xi}}(t))\leq d(o,p)$ for all $t\geq 0$. Thus
\begin{align*}
(x|y)_{o,\xi}&=\lim_{t\to\infty}t-\frac{1}{2}d(\gamma_{ox_{\xi}}(t),\gamma_{oy_{\iota\xi}}(t)) \\
& \geq \lim_{t\to\infty}t-\frac{1}{2}d(\gamma_{px_{\xi}}(t),\gamma_{py_{\iota\xi}}(t))-d(o,p)=\infty.
\end{align*}

Let $(x,y)\in \Ftm\times\Ftim$ be such that $x\op y$. Then by the above proposition  $(x|y)_{o,\xi} =\frac{1}{2}b_{x_{\xi}}(o,\rho_{c_x,c_y}(o))\leq d(o,\rho_{c_x,c_y}(o))$, i.e. $(x|y)_{o,\xi}< \infty$.
\end{proof}

The above corollary implies that $\Af$ is the maximal domain of definition for $\cro_{o,\xi}$. As mentioned, in case of a symmetric space $X$ is the set $\Af$ an open and dense subset of $(\Ftm(X)\times \Ftim(X))^2$, i.e. the cross ratio is generically defined.

\bp\label{proposition basepoint change} Let $o,\hat{o}\in M$, $(x,y)\in \Ftm\times \Ftim$ and $\xi\in\tm$. Then 
$(x|y)_{o,\xi}=(x|y)_{\hat{o},\xi}+\frac{1}{2}b_{x_{\xi}}(o,\hat{o})+\frac{1}{2}b_{y_{\iota\xi}}(o,\hat{o}).$
\ep

\begin{proof}
If $x\nop y$, then by the above corollary $(x|y)_{o,\xi}=\infty=(x|y)_{\hat{o},\xi}$. 

If $x\op y$, let $\rx, \ry$ be any horospherical retracts as in Proposition \ref{prop Gromov product as Busemann function}. Then 
\begin{align*}
b_{x_{\xi}}(o,\ry(o))=b_{x_{\xi}}(o,\hat{o})+b_{x_{\xi}}(\hat{o},\ry(\hat{o}))+b_{x_{\xi}}(\ry(\hat{o}),\ry(o)).
\end{align*}  

By construction $\ry(o) ,\ry(\hat{o})\in P(x,y)$. Moreover $x,y$ are opposite and hence by Lemma \ref{lem busemann in parallel set} and equation \eqref{eq horospherical retract and busemann} 
\begin{align*}
b_{x_{\xi}}(\ry(\hat{o}),\ry(o))=-b_{y_{\iota\xi}}(\ry(\hat{o}),\ry(o))=-b_{y_{\iota\xi}}(\hat{o},o)=b_{y_{\iota\xi}}(o,\hat{o}).
\end{align*}

Together with Proposition \ref{prop Gromov product as Busemann function} the claim follows.
\end{proof}

\bp\label{prop cr base point independence}
Let $o,\hat{o}\in M$. Then $\cro_{o,\xi}(x_1,y_1,x_2,y_2)=\cro_{\hat{o},\xi}(x_1,y_1,x_2,y_2)$ for all $(x_1,y_1,x_2,y_2)\in \Af$.
\ep 

\begin{proof}
Plugging in the above proposition in the definitions of $\cro_{o,\xi}$ and $\cro_{\hat{o},\xi}$ yields directly the result.
\end{proof}

\bd\label{def cross ratio}
Given $(x_1,y_1,x_2,y_2)\in \Af$, we define the \emph{cross ratio with respect to $\xi\in \sm$} to be $\cro_{\xi}(x_1,y_1,x_2,y_2)=\cro_{o,\xi}(x_1,y_1,x_2,y_2)$ for some $o\in M$.
\ed

\begin{Example}\label{exam h2 cross h2}(see also \cite{Kim}) Consider the symmetric space $X=\Ht\times\Ht$, where $\Ht$ is the hyperbolic plane. The ideal boundary $\di(\Ht\times\Ht)$ can be identified with $S^1\times S^1\times [0,\frac{\pi}{2}]$ - this is in such a way that the unit-speed geodesic ray from a base-point $(o_1,o_2)\in\Ht\times\Ht$ to the point in $(x_1,x_2,\alpha)\in S^1\times S^1\times [0,\frac{\pi}{2}]\cong\di(\Ht\times\Ht)$ is given by $(\gamma_{o_1x_1}(\cos(\alpha)t), \gamma_{o_2x_2}(\sin(\alpha)t))$.

The types are exactly determined by the angle $\alpha$ and the opposition involution equals the identity. In particular every type is self opposite.

Fix $o=(o_1,o_2)\in\Ht\times\Ht$ and $x=(x_1,x_2,\alpha), y=(y_1,y_2,\alpha)\in \di(\Ht\times\Ht)$ and set $\gamma_1:=\gamma_{o_1 x_1}, \hat{\gamma}_1:=\gamma_{o_1 y_1}, \gamma_2:=\gamma_{o_2 x_2}$ and $\hat{\gamma}_2:=\gamma_{o_2 y_2}$. Then
\begin{align*}
(x|y)_{o,\alpha}=\lim\limits_{t\to\infty} t-\frac{1}{2}\sqrt{|\gamma_1(\cos(\alpha) t)\hat{\gamma}_1(\cos(\alpha) t)|^2+|\gamma_2(\sin(\alpha) t)\hat{\gamma}_2(\sin(\alpha) t)|^2}.
\end{align*}
Using $\lim_{t\to\infty}|\gamma_1(\cos(\alpha) t)\hat{\gamma}_1(\cos(\alpha) t)| -2\cos(\alpha)t=-2(x_1|y_1)_{o_1}$, if $\alpha\neq \frac{\pi}{2}$
\begin{align*}
& (x|y)_{o,\alpha} =\lim\limits_{t\to\infty} t-\sqrt{(-(x_1|y_1)_{o_1}+\cos(\alpha) t)^2+(-(x_2|y_2)_{o_2}+\sin(\alpha) t)^2}\\
= & \lim\limits_{t\to\infty} t-\sqrt{t^2-2 t(\cos(\alpha)(x_1|y_1)_{o_1}+\sin(\alpha)(x_2|y_2)_{o_2})+(x_1|y_1)^2_{o_1}+(x_2|y_2)^2_{o_2}}.
\end{align*}
We substitute $t=s^{-1}$. Then a Taylor expansion for the root expression at $s=0$ yields that 
\begin{align*}
(x|y)_{o,\alpha}&=\lim\limits_{s\to 0} \frac{1}{s}(1-(1- s(\cos(\alpha)(x_1|y_1)_{o_1}+\sin(\alpha)(x_2|y_2)_{o_2}) + o(s))\\
&=\cos(\alpha)(x_1|y_1)_{o_1}+\sin(\alpha)(x_2|y_2)_{o_2}.
\end{align*}

Therefore $\cro_{\alpha} =\cos(\alpha) \log |\cro_{\Ht}|+\sin(\alpha)\log|\cro_{ \Ht}|$, where $\cro_{\Ht}$ is the usual \emph{multiplicative} cross ratio on $\di\Ht$.
\end{Example}

\bl\label{lem Gromov and cro product continuous} Let $X$ be a symmetric space. Then for every $o\in X$ the Gromov product $(\cdot|\cdot)_{o,\xi}:\Ftm(X)\times\Ftim(X) \to [0,\infty]$ is continuous. In particular also $\cro_{\xi}$ is continuous.
\el

\begin{proof} 
Since $\Ftm(X),\Ftim(X)$ are manifolds it is enough to consider sequential continuity. Therefore let $(x,y)\in \Ftm(X)\times\Ftim(X)$ and let $x_i\to x$ and $y_i\to y$. 

If $x\nop y$, we have $(x|y)_{o,\xi}=\infty$. 
We set $(x|y)_{o,\xi}(t):=(\gamma_{ox_{\xi}}(t)|\gamma_{oy_{\iota\xi}}(t))_o$ with Gromov product on the right hand side the usual Gromov product on the metric space $(X,d)$. As $X$ is non-positively curved, the function $t\mapsto (x|y)_{o,\xi}(t)$ is monotone increasing.
Let $C>0$ be given. Then there is $t_C\in \R_{+}$ such that $(x|y)_{o,\xi}(t_C)\geq C+2$. 
Since the topology on $\Ftm(X)$ is induced by the cone topology, we have that $(x_i)_{\xi}\to x_{\xi}$ in the cone topology and similarly for $y_i$ and $y$.
Hence we find $L\in\N$ such that $d(\gamma_{o(x_i)_{\xi}}(t_C),\gamma_{ox_{\xi}}(t_C))<1$ and $d(\gamma_{o(y_i)_{\iota\xi}}(t_{C}),\gamma_{oy_{\iota\xi}}(t_{C}))<1$ for all $i\geq L$.
Hence by the triangle inequality $(x_i|y_j)_{o,\xi}(t_C)>(x|y)_{o,\xi}(t_C)-2> C$ for all $i,j\geq L$. As $C$ was arbitrary, this yields $\lim_{i,j\to\infty}(x_i|y_j)_{o,\xi}=\infty$ - which proves continuity for $x\nop y$.

Assume $x\op y$. Let $K=\st_G(o)$. We know that $K$ acts transitively on $\Ftm(X)$ and we have a $K$-equivariant and homeomorphic identification $\Ftm(X)\simeq K\slash K_x$. Therefore $x_i\to x$ implies that we find $k_i\in K$ such that $k_i x_i=x$ and $k_i\to e\in G$. Now, $x\op y$ and opposition being an open condition, together with $y_i\to y$ and $k_i\to e$, imply that there exists $L\in \N$ such that $k_iy_j\op x$ for all $i,j\geq L$. 
Thus there exists a unique $n_{ij}\in N_x$ such that $n_{ij}k_iy_j=y$ for $i,j\geq L$. From $k_i\to e$ and $y_j\to y$ it follows $n_{ij}\to e\in G$ for $i,j\to\infty$. We set $g_{ij}:=n_{ij}k_i$ and by construction $g_{ij}\to e$, $g_{ij}x_i=x$, $g_{ij}y_j=y$. Hence $(x_i|y_j)_{o,\xi}=(x|y)_{g_{ij}o,\xi}$. Proposition \ref{proposition basepoint change} and $g_{ij}\to e$ yield that $(x_i|y_j)_{o,\xi}\to (x|y)_{o,\xi}$.
\end{proof}

\bl\label{lem gromov product continuous in types}
Let $(x,y)\in \Ftm\times\Ftim$ and $x\op y$. Moreover, let $\xi_i\in\tm$ be a sequence with $\xi_i\to \xi_0\in\tm$. Then $(x|y)_{o,\xi_i}\to (x|y)_{o,\xi_0}$. In particular, $\cro_{\xi_i}(x,y,z,w)\to \cro_{\xi_0}(x,y,z,w)$ for all $(x,y,z,w)\in\Atop$.
\el

\begin{proof}
Let $c_x,c_y\in\Ftsm$ such that $c_x\op c_y$, $x$ is a face of $c_x$ and $y$ is a face of $c_y$. Then Proposition \ref{prop Gromov product as Busemann function} and equation \eqref{eq horospherical retract and busemann} imply $(x|y)_{o,\xi}=\frac{1}{2} b_{x_{\xi}}(\rho_{c_x,c_y}(o),\rho_{c_y,c_x}(o))$ for all $\xi\in\tm$. Denote $p_x:=\rho_{c_x,c_y}(o)$, $p_y:=\rho_{c_y,c_x}(o)$ and by $A_{xy}$ the unique affine apartment with $c_x,c_y\subset \di A_{xy}$. 

Every affine apartment can be isometrically identified with $\R^r$ where $r$ is the rank of $M$. We identify $A_{xy}$ with $\R^r$ such that $0\simeq p_x$. Let $v_{\xi}\in A_{xy}\simeq \R^r$ be of norm one and such that the line from $0$ through $v_{\xi}$ is the geodesic ray in $A_{xy}$ from $p_x$ to $x_{\xi}$. Then Euclidean geometry yields that $b_{x_{\xi}}(p_x,p_y)=\langle v_{\xi}, p_y\rangle$. In particular, we get 
\begin{align}\label{eq gromov product as busemann in apartment}
(x|y)_{o,\xi_i}=\frac{1}{2} \langle v_{\xi_i}, p_y\rangle. 
\end{align}
Moreover $\xi_i\to\xi_0$ implies that $v_{\xi_i}\to v_{\xi_0}$ and hence the claim follows.
\end{proof}

The assumption of opposition in the above lemma is needed, since there are $(x,y)\in \Ftm\times\Ftim$ with $x\nop y$ but there are faces $x_0$ of $x$ and $y_0$ of $y$ with $x_0\op y_0$. Then if $\xi_i\in\inn(\tm)$ converge to $\xi_0$ such that $\xi_0\in\inn(\tau_0)$ and $\tau_0$ is the type of $x_0$, we get $(x|y)_{o,\xi_i}=\infty \nrightarrow (x_0|y_0)_{o,\xi_0}$ (as the latter is finite).

We remind that any isometry $F:M_1\to M_2$ induces a building isomorphism $F_{\infty}:\Di M_1\to \Di M_2$ together with a type map $F_{\sm}:\sm_1\to\sm_2$ with the property that $F(\Ftm(M_1))=\mathrm{Flag}_{F_{\sm}(\tm)}(M_2)$. 

\bp\label{prop isometries are Moebius maps}
Let $F:M_1\to M_2$ be an isometry between either symmetric spaces or thick Euclidean buildings, $F_{\infty}:\Di M_1\to \Di M_2$ the induced building isomorphism and $\xi\in\sm_1$. Then
\begin{align*}
\cro_{\xi_1}(x_1,y_1,x_2,y_2)=\cro_{F_{\sm}(\xi_1)}(F_{\infty}( x_1),F_{\infty}(y_1),F_{\infty}(x_2),F_{\infty}(y_2))
\end{align*}
for all $(x_1,y_1,x_2,y_2)\in \mathcal{A}_{\xi_1}$. Equivalently, $\cro_{\xi_1}=F_{\infty}^{*}\cro_{F_{\sm}(\xi_1)}$ with $F_{\infty}^{*}$ denoting the pullback under $F_{\infty}$.
\ep

\begin{proof}
Let $\xi_1\in \tm$ and $(x,y)\in \Ftm(M_1)\times\Ftim(M_1)$. Since the Gromov product $(\cdot|\cdot)_{o,\xi_1}$ is defined in terms of a limit of distances involving unit speed geodesics and isometries leave those invariant, it follows that $(x|y)_{o,\xi_1}=(F_{\infty}(x)|F_{\infty}(y))_{F(o),F_{\sm}(\xi_1)}$.  Hence Corollary \ref{corollary-not-opposite-equivalent-distance-zero} implies that if $(x_1,y_1,x_2,y_2)\in \mathcal{A}_{\xi_1}$, then $(F_{\infty}(x_1),F_{\infty}(y_1),F_{\infty}(x_2),F_{\infty}(y_2))\in \mathcal{A}_{F_{\sm}(\xi_1)}$. Finally, $\cro_{\xi_1}=\cro_{o,\xi_1}=F_{\infty}^{*}\cro_{F(o),F_{\sm}(\xi_1)}=F_{\infty}^{*}\cro_{F_{\sm}(\xi_1)}$ by Proposition \ref{prop cr base point independence}.
\end{proof}

\bc
Let $g\in Iso(M)$ and $\xi_0$ be the center of gravity of $\sm$ with respect to the angular metric. Then $\cro_{\xi_0}=g^{*}\cro_{\xi_0}$. In case of a symmetric space $X$ and $g\in G$ we have $\cro_{\xi,X}=g^{*}\cro_{\xi,X}$ for all $\xi\in\sm$.
\ec

\begin{proof}
For the center of gravity $\xi_0\in\sm$ we have $g_{\sm}(\xi_0)=\xi_0$ for all $g\in Iso(M)$, as $g_{\sm}:\sm\to\sm$ is an isometry with respect to the angular metric. Then the first claim follows. In case of a symmetric space and $g\in G$, we know $g_{\sm}=id_{\sm}$, which implies the second claim.
\end{proof}

\begin{Example}\label{ex Sl n}
We want to determine the Gromov products and cross ratios of the symmetric spaces $X(n):=\text{SL}(n,\R)\slash \text{SO}(n,\R)$. For a deeper description of the symmetric space $X(n)$ see \cite{Eberlein}.\medskip


The ideal boundary $\di X(n)$ can be identified with eigenvalue flag pairs $(\lambda,F)$, where $F=(V_1,\ldots,V_l)$ is a flag in $\R^n$, i.e. the $V_i$ are subspaces of $\R^n$ with $V_i\varsubsetneq V_{i+1}$, $V_l=\R^n$, and $\lambda=(\lambda_1,\ldots,\lambda_l)\in \mathbb{R}^l$ such that $\lambda_i>\lambda_{i+1}$, $\sum_{i=1}^{l} m_i \lambda_i=0$ for $m_i=\dim V_{i} - \dim V_{i-1}$ and $\sum_{i=1}^l m_i \lambda^2_i=1$. In particular, $2\leq l \leq n$.
The action of $g\in \text{SL}(n,\R)$ on an eigenvalue flag pair is given by $g\cdot (\lambda,F)=(\lambda,g\cdot F)$, where $g\cdot (V_1,\ldots,V_l)=(g\cdot V_1,\ldots,g\cdot V_l)$ and $F=(V_1,\ldots,V_l)$.

The "eigenvalues" $\lambda$ in the eigenvalue flag pairs $(\lambda,F)$ determine the type of any point in the ideal boundary. More precisely, the set of paris $(\lambda_1,\ldots,\lambda_l),(m_1,\ldots,m_l)$, $\lambda_i\in\R,m_i\in\N\backslash\{0\}$ with $\lambda_i>\lambda_{i+1}$, $\sum_{i=1}^{l} m_i \lambda_i=0$, $\sum_{i=1}^l m_i \lambda^2_i=1$ and $\sum_{i=1}^l m_i=n$ parametrize the Weyl chamber $\sm$. We have that $\lambda=(\lambda_1,\ldots,\lambda_l)$ is in the interior of the chamber if and only if $l=n$.

Faces of $\sm$ can be charaterized in the following way: Two pairs as above $(\lambda_1,\ldots,\lambda_l),(m_1,\ldots,m_l)$ and $(\lambda'_1,\ldots,\lambda'_l),(m'_1,\ldots,m'_l)$ are in the interior of the same face if and only if $m_i=m_i'$ for all $i=1,\ldots,l$. In particular we can identify the set of faces of $\sm$ with $\{(m_1,\ldots,m_l)\in\N^{l} | l \geq 2, m_i\neq 0, \sum_{i=1}^{l} m_i = n \}$. For $\tm \simeq(m_1,\ldots,m_l)$ we have $\Ftm=\{(V_1,\ldots,V_l)| V_i\varsubsetneq V_{i+1}, \dim V_i - \dim V_{i-1}=m_i\}$.
The action of the opposition involution is given by $\iota(\lambda_1,\ldots,\lambda_l)=(-\lambda_l,\ldots,-\lambda_1)$ and $\iota (m_1,\ldots,m_l)=(m_l,\ldots,m_1)$. Hence, if $V=(V_1,\ldots,V_l)\in\Ftm$ and $W=(W_1,\ldots,W_l)\in\Ftim$, then $\dim V_i + \dim W_{l-i}=n$. In this situation $V\op W\Longleftrightarrow V_i \oplus W_{l-i}=\R^n$ for all $i=1,\ldots,l-1$.

Let $V=(V_1,\ldots,V_l),Y=(Y_1,\ldots,Y_l)\in\Ftm$ and $W=(W_1,\ldots,W_l),Z=(Z_1,\ldots,Z_l)\in\Ftim$ such that $V,Y\op W,Z$. Let $i_j=\dim V_j$. Then fix a basis $(v_1,\ldots v_n)$ such that $V_j=\text{span}\lbrace v_1,\ldots, v_{i_j} \rbrace$. In the same way we fix basis $(w_1,\ldots w_n), (y_1,\ldots y_n)$ and  $(z_1,\ldots z_n)$ for $W,Y,Z$, respectively. 
Additionally, fix an identification $\wedge^n \R^n\cong \R$. We set $V_j\wedge W_{l-j}:=v_1\wedge \ldots\wedge v_{i_j}\wedge w_1 \wedge \ldots\wedge w_{n-i_j}$ (we have $W_{l-j}=\text{span}\lbrace w_1,\ldots, w_{n-i_j}\rbrace$) and in the same way for the other flags. Then the term
$
(V_j\wedge W_{l-j})(Y_j\wedge Z_{l-j})(V_j\wedge Z_{l-j})^{-1}(Y_j\wedge W_{l-j})^{-1}
$
can be shown to be independent of all choices for all $j=1,\ldots, l-1$ - compare e.g. \cite{Martone-Zhang}.

Let $V,W,Y,Z$ be as before and $\lambda=(\lambda_1,\ldots,\lambda_l)$ a type with $\lambda\in\inn(\tm)$. Then
\begin{align*}
\cro_{\lambda}(V,W,Y,Z)= n \sum_{j=1}^{l-1}(\lambda_j - \lambda_{j+1})\log(|\frac{V_j\wedge W_{l-j}}{V_j\wedge Z_{l-j}} \frac{Y_j\wedge Z_{l-j}}{Y_j\wedge W_{l-j}}|),
\end{align*}
using the above conventions - see the appendix for a proof. We remark that some specific of those cross ratios are known already and have been used for analyzing  Hitchin representations and more general Anosov representations (see e.g. \cite{Labourie-IHES}, \cite{Martone-Zhang}).
\end{Example}

Let $M=M_1\times \ldots\times M_k$ be a product of either symmetric spaces or Euclidean buildings. Then the building at infinity $\Di M$ is the spherical join of the buildings $\Di M_i$ \cite[Sc.4.3]{Kleiner-Leeb}. In particular, the Weyl chamber $\sm$ decomposes as a spherical join $\sm=\sm_1\circ \ldots\circ \sm_k$. Hence we get a surjective map
\begin{align}\label{eqn types in product spaces}
\pi:\sm_1\times\ldots\times \sm_k\times S_k^{+}\to \sm,
\end{align}
where $S_k^{+}:=\lbrace \mu=(\mu_1,\ldots,\mu_k) \in [0,1]^k\;|\; \Sigma_1^k \mu_i^2 =1\rbrace$. We remark that $\pi$ is in general not injective, since it is independent of the exact choice of the type $\xi_i\in \sm_i$ if $\mu_i=0$.

Let $\xi=\pi(\xi_1,\ldots,\xi_k,\mu)$ with $\mu=(\mu_1,\ldots,\mu_k)\in S_k^{+}$ and let $x=(x_1,\ldots,x_k)\in \Ftm(M)\simeq \mathrm{Flag}_{\tau_1}(M_1)\times\ldots\times \mathrm{Flag}_{\tau_k}(M_k)$\footnote{Actually we would have a spherical join instead of the product. However, we can naturally identify a simplex in a join with the product of the simplices in the different factors - and that is what we do here for simplicity.}
 such that $\xi\in\inn(\tm)$ and $\xi_i\in\inn(\tm_i)$. For simplicity we assume $\mu_i\neq 0$ for all $1\leq i\leq k$ - if some $\mu_i=0$ essentially the same formula holds, but the factor $\mathrm{Flag}_{\tau_i}(M_i)$ is not apparent in the decomposition of $\Ftm(M)$. 

We remark that the unit-speed geodesic from some point $(o_1,\ldots,o_k)\in M$ to $x_{\xi}$ is of the form $(\gamma_{o_1x_{\xi_1}}(\mu_1 t),\ldots,\gamma_{o_k x_{\xi_k}}(\mu_k t))$, where $\gamma_{o_i x_{\xi_i}}$ denote the unit speed geodesics in the factors $M_i$ joining $o_i$ to $(x_i)_{\xi_i}$ - cp. also Example \ref{exam h2 cross h2}.

Let $y=(y_1,\ldots,y_k)\in  \Ftim(M)\simeq \mathrm{Flag}_{\iota\tau_1}(M_1)\times\ldots\times \mathrm{Flag}_{\iota\tau_k}(M_k)$ and be $x$ and $\xi$ as above. Then similar calculations as in Example \ref{exam h2 cross h2}, yield that
\begin{align*}
(x|y)_{(o_1,\ldots, o_k),\pi(\xi_1,\ldots,\xi_k,\mu)}=\mu_1 (x_1|y_1)_{o_1,\xi_1}+\ldots+\mu_k (x_k|y_k)_{o_k,\xi_k}.
\end{align*}

\bp\label{pro cross ratio as product}
Notations as before. Moreover, let $z\in \Ftm(M)$ and $w\in\Ftim(M)$. Then
\begin{align*}
\cro_{\pi(\xi_1,\ldots,\xi_k,\mu)}(x,y,z,w)
=  \mu_1\cro_{\xi_1}(x_1,y_1,z_1,w_1)+\ldots+ \mu_k \cro_{\xi_k}(x_k,y_k,z_k,w_k)
\end{align*}
for $(x,y,z,w)\in \mathcal{A}_{\pi(\xi_1,\ldots,\xi_k,\mu)}$.
\ep

\section{Vector valued cross ratios}

So far, we have constructed families of cross ratios on subsets of the spaces $(\Ftm\times\Ftim)^2$ which are parametrized by $\xi\in\inn(\tm)$. In this section we show that such a family gives rise to a single vector valued cross ratio containing all the information of the family. The vector valued cross ratio has the same symmetries as the usual cross ratios (cp. equations~\eqref{eq symmetries of cross ratios}) justifying the name cross ratio.

We remind that $\sm=\apo$; hence every type can be viewed as vector in $\am$ of norm one.

\bl\label{lem gromov prod as sum of others}
Let $\tm$ be a face of $\sm$ and $\xi_0,\xi_1,\ldots, \xi_j\in\tm$ such that there exist $a_i\in\R$ with $\xi_0= \sum_{i=1}^j a_i \xi_i$. Then for $(x,y)\in\Ftm\times \Ftim$ with $x\op y$ we have 
$(x|y)_{o,\xi_0} =\sum_{i=1}^j a_i (x|y)_{o,\xi_i}$.

In particular it follows that $\cro_{\xi_0}(x,y,z,w) =\sum_{i=1}^j a_i \cro_{\xi_i}(x,y,z,w)$ for all $(x,y,z,w)\in\Atop$
\el

\begin{proof}
Let $c_x,c_y\in\Ftsm$ such that $c_x\op c_y$, $x$ is a face of $c_x$ and $y$ is a face of $c_y$. We recall the notation of the proof of Lemma \ref{lem gromov product continuous in types}: We denote $p_x:=\rho_{c_x,c_y}(o)$, $p_y:=\rho_{c_y,c_x}(o)$ and by $A_{xy}$ the unique apartment with $c_x,c_y\subset\di A_{xy}$. Moreover, let $A_{xy}\simeq \R^r$ such that $p_x \simeq 0$, in particular $A_{xy}$ inherits a inner product. Let $v_{\xi}\in A_{xy}\simeq \R^r$ be of norm one and such that the line from $p_x\simeq 0$ through $v_{\xi}$ is the geodesic ray in $A_{xy}$ from $p_x$ to $x_{\xi}$. Then we know from equation \eqref{eq gromov product as busemann in apartment} that $(x|y)_{o,\xi_i}=\frac{1}{2} \langle v_{\xi_i}, p_y\rangle$.

By the definition of the $v_{\xi_i}$ it is immediate that $v_{\xi_0}= \sum_{i=1}^j a_i v_{\xi_i}$, where we have the addition inherited to $A_{xy}$ under the identification with $\R^r$ such that $p_x\simeq 0$. Hence
\[
(x|y)_{o,\xi_0}= \frac{1}{2}\langle v_{\xi_0}, p_y\rangle=\sum_{i=1}^j \frac{1}{2}a_i \langle v_{\xi_i}, p_y\rangle=\sum_{i=1}^j a_i (x|y)_{o,\xi_i}. \qedhere
\]
\end{proof}

Let $\xi_1,\ldots, \xi_r\in\am$ be the corners of $\sm=\apo$. Then every subset $J\subset \lbrace 1, \ldots,r\rbrace$ defines a simplex in $\sm$, i.e. a face $\tm$ of $\sm$. In the same way every simplex $\tm\subset\sm$ gives a subset $J_{\tm} \subset \lbrace 1, \ldots,r\rbrace$. 

Given a simplex $\tm$ we recall that $\at=\text{span}_{j\in J_{\tm}} \xi_j\subset \am$. Moreover, we define $\alpha_j^{\tm}\in \at$ for $j\in J_{\tm}$ by $\langle \alpha_j^{\tm}, \xi_i\rangle=\delta_{ij}$ for all $i\in J_{\tm}$ - this yields well defined vectors, as the $\xi_i$ with $i\in J_{\tm}$ form a basis of $\at$. We remind that $\am$ was naturally equipped with an inner product.

The $\xi_j$ correspond to normalized fundamental weights of the root system and the $\alpha^{\sm}_j$ to possibly rescaled roots.

\bd
Let $\tm$ be a face of $\sm$ and $J_{\tm}$, $\alpha_j^{\tm}$ as above. Then we define a (vector valued) cross ratio $\cro_{\tm}:\At\to \at\cup\lbrace \pm\infty\rbrace$ by 
\begin{align*}
\cro_{\tm}(x,y,z,w):=\sum_{i\in J_{\tm}} \cro_{\xi_i}(x,y,z,w) \alpha^{\tm}_i.
\end{align*}
Here we set $\cro_{\tm}(x,y,z,w):=-\infty$ if $x\nop y$ or $z\nop w$ and $\cro_{\tm}(x,y,z,w):=\infty$ if $x\nop w$ or $z\nop y$. 
\ed

It is straight forward to see that $\cro_{\tm}$ has the same symmetries as in equations \eqref{eq symmetries of cross ratios}, where the addition is now in the vector space $\at$.

The vector valued cross ratio contains the full information of the collection of cross ratios form the previous section:

\bl\label{lem vector valued cro gives usual by inner product} Let $\xi\in\inn(\tm)$. Then $\langle\cro_{\tm}(x,y,z,w),\xi\rangle=\cro_{\xi}(x,y,z,w)$ for $(x,y,z,w)\in \Atop$ and $\cro_{\tm}(x,y,z,w)=\pm\infty=\cro_{\xi}(x,y,z,w)$ for $(x,y,z,w)\in \At\backslash \Atop$.
\el

\begin{proof}
If $(x,y,z,w)\in \At\backslash \Atop$, then the equality is immediate. Hence assume $(x,y,z,w)\in \Atop$. Then
\begin{align*}
\langle\cro_{\tm}(x,y,z,w),\xi\rangle=\sum_{i\in J_{\tm}} \cro_{\xi_i}(x,y,z,w)\langle \alpha^{\tm}_i,\xi\rangle.
\end{align*}
Since $\langle \alpha_j^{\tm}, \xi_i\rangle=\delta_{ij}$ for all $i\in J_{\tm}$, we derive that $\langle\sum_{i\in J_{\tm}} \langle \alpha^{\tm}_i,\xi\rangle\xi_i,\alpha^{\tm}_j\rangle=\langle\xi, \alpha^{\tm}_j\rangle$ for all in $j\in J_{\tm}$. Moreover, it is immediate that the $\alpha^{\tm}_j$ form a base of $\at$. Thus we get that $\sum_{i\in J_{\tm}} \langle \alpha^{\tm}_i,\xi\rangle\xi_i=\xi$. Therefore Lemma \ref{lem gromov prod as sum of others} implies $\sum_{i\in J_{\tm}} \langle \alpha^{\tm}_i,\xi\rangle \cro_{\xi_i}(x,y,z,w)=\cro_{\xi}(x,y,z,w)$.
\end{proof}

The above lemma also holds for $\xi\in\partial \tm$ as long as $(x,y,z,w)\in\Atop$, but does not hold for general $(x,y,z,w)\in \At$ - in this case $\cro_{\xi}(x,y,z,w)$ might be finite while $\cro_{\tm}(x,y,z,w)$ is not (compare the discussion just after Lemma \ref{lem gromov product continuous in types}).

The following corollary captures the topological properties of $\cro_{\tm}$ in case of symmetric spaces. It is an immediate consequence of the lemma above and Lemma \ref{lem Gromov and cro product continuous}.

\bc
Let $X$ be a symmetric space. The map $\cro_{\tm}$ restricted to $\Atop$ is continuous and for all $\xi\in \inn(\tm)$ the map $\langle\cro_{\tm}(\cdot),\xi\rangle:\At\to \R\cup\lbrace \pm\infty\rbrace$ is continuous.
\ec

Let $\pi_{\tm}:\am\to\at$ be the orthogonal projection. Then it is straight forward to show that $\pi_{\tm}(\alpha^{\sm}_i)=\alpha^{\tm}_i$ for all $i\in J_{\tm}$ and $\pi_{\tm}(\alpha^{\sm}_j)=0$ for all $j\notin J_{\tm}$. Then we can derive that $\cro_{\tm}(x,y,z,w)= \pi_{\tm}(\cro_{\sm}(x,y,z,w))$ for all $(x,y,z,w)\in\Asop$.

\subsection*{Translation vectors and periods}

We assume for this section that $\tm$ is \textbf{self-opposite}, i.e. $\tm=\iota\tm$. Moreover denote by $Iso_{e}(M)$ the subgroup of $Iso(M)$ such that $g_{\sm}=\mathrm{id}$ for all $g\in Iso_{e}(M)$ - in particular $G=Iso_{e}(X)$ for a symmetric space $X$. Let $g\in Iso_{e}(M)$ such that $g$ stabilizes two points $g^{\pm}\in \Ftm$ with $g^{-}\op g^{+}$. Since $g$ is an isometry, it maps every geodesic connecting points of the interior of $g^{-}$ and $g^{+}$ to another geodesic connecting the same points. In particular $g$ stabilizes $P(g^{-},g^{+})$ set-wise.

In the preliminaries we have seen that $P(g^{-},g^{+})$ splits as a product $\at\times CS(g^{-},g^{+})$ such that $g^{\pm}$ are identified with the positive and negative, respectively, maximal dimensional simplices in $\at$, i.e. $g^{+}\simeq\di \at^{+}$ where $\at^{+}:=\at\cap \am^{+}$. 
Note that $g$ descends to an isometry $g_{\at}$ of $\at$. Since $\at$ is Euclidean and $g_{\at}$ stabilizes each boundary point of $\at$, $g_{\at}$ acts as a translation on $\at$. More precisely, there exists a \emph{translation vector} $\ell_{g}^{\tm}\in\at$ such that $g_{\at}(p)=p+\ell_{g}^{\tm}$ for all $p\in\at$. 

\bp\label{pro periods of isometries} Let $g\in Iso_e(M)$ such that $g^{\pm}\in \Ftm$ with $g^{-}\op g^{+}$ are stabilized by $g$. Let $\ell_{g}^{\tm}$ denote the translation vector along the first factor of $P(g^{-},g^{+})\simeq \at\times CS(g^{-},g^{+})$. Then
$\cro_{\tm}(g^{-},g\cdot x,g^{+},x)= \frac{1}{2}(\ell_{g}^{\tm}+\iota \ell_{g}^{\tm})$,
for any $x\in\Ftm$ with $x\op g^{\pm}$.
\ep 

\begin{proof}
We remark that $\cro_{\tm}(g^{-},g\cdot x,g^{+},x)$ is independent of the choice of $x\op g^{\pm}$; this follows from the symmetries of $\cro_{\tm}$ together with Proposition \ref{prop isometries are Moebius maps}. Therefore, we fix one $x\in\Ftm$ with $x\op g^{\pm}$.

Let $o\in P(g^{-},g^{+})$ and $\xi_i$ with $i\in J_{\tm}$ be the corners of $\tm$. By assumption $x\op g^{\pm}$ and hence $g\cdot x\op g^{\pm}$. Then Proposition \ref{proposition basepoint change} yields
\begin{align*}
(g^{\pm}|g\cdot x)_{o,\xi_i}=(g^{\pm}| x)_{g^{-1}\cdot o,\xi_i}= (g^{\pm}| x)_{ o,\xi_i}+\frac{1}{2} b_{g^{\pm}_{\xi_i}}(g^{-1}\cdot o,o)+\frac{1}{2} b_{x_{\iota\xi_i}}(g^{-1}\cdot o,o).
\end{align*}
Moreover, we have $b_{g^{\pm}_{\xi_i}}(g^{-1}\cdot o,o)= b_{g^{\pm}_{\xi_i}}(o,g\cdot o)$. If we plug this in the definition of $\cro_{\xi_i}$, several terms cancel and we are left with $\cro_{\xi_i}(g^{-},g\cdot x,g^{+},x)= \frac{1}{2}b_{g_{\xi_i}^{+}}(o,g\cdot o)- \frac{1}{2}b_{g_{\xi_i}^{-}}(o,g\cdot o)$.
Since $o,g\cdot o\in P(g^{-},g^{+})$ and $ g_{\iota\xi_i}^{+}\in g^{+}$ is the point opposite to $g_{\xi_i}^{-}\in g^{-}$, Lemma \ref{lem busemann in parallel set} implies $b_{g_{\xi_i}^{-}}(o,g\cdot o)=-b_{g_{\iota\xi_i}^{+}}(o,g\cdot o)$. In particular $\cro_{\xi_i}(g^{-},g\cdot x,g^{+},x)= \frac{1}{2}b_{g_{\xi_i}^{+}}(o,g\cdot o)+ \frac{1}{2}b_{g_{\iota\xi_i}^{+}}(o,g\cdot o)$.

Since $o$ was arbitrary in $P(g^{-},g^{+})$ we can assume that its first coordinate under the identification $P(g^{-},g^{+})\simeq\at\times CS(g^{-},g^{+})$ is $0\in\at$. Moreover, we can use Lemma \ref{lem horospheres in parallel set} to see that only the first factor matters for the Busemann functions $b_{g_{\xi_i}},b_{g_{\iota\xi_i}}$. As $g$ acts as a translation on $\at$, we have that $g\cdot 0= \ell_{g}^{\tm}$. Hence $b_{g_{\xi_i}^{+}}(o,g\cdot o)=\langle \xi_i, \ell_{g}^{\tm}\rangle$ (cp. the arguments around equation \eqref{eq gromov product as busemann in apartment}). 
By assumption $\tm=\iota\tm$, hence $\iota$ restricts to an isometry $\iota:\at\to\at$. Together with $\iota^2=id$, this yields $\langle \iota \xi_i,\ell_{g}^{\tm}\rangle=\langle \xi_i,\iota\ell_{g}^{\tm}\rangle$.
Altogether we derive
\begin{align*}
\cro_{\tm}(g^{-},g\cdot x,g^{+},x)= \sum_{i\in J_{\tm}}\frac{1}{2}(\langle \xi_i,\ell_{g}^{\tm}\rangle+ \langle \xi_i,\iota\ell_{g}^{\tm}\rangle) \alpha^{\tm}_i.
\end{align*}
It is immediate that $\langle \cro_{\tm}(g^{-},g\cdot x,g^{+},x),\xi_i\rangle=\frac{1}{2}(\langle \xi_i,\ell_{g}^{\tm}\rangle+ \langle \xi_i,\iota\ell_{g}^{\tm}\rangle)$ for all $i\in J_{\tm}$. Since the $\xi_i$ with $i\in J_{\tm}$ form a basis of $\tm$, it follows that $\cro_{\tm}(g^{-},g\cdot x,g^{+},x)= \frac{1}{2}(\ell_{g}^{\tm}+\iota\ell_{g}^{\tm})$.
\end{proof}

Let $g\in Iso_e (M)$ be as before. Then the term $\cro_{\tm}(g^{-},g\cdot x,g^{+},x)$ is also called \emph{period} - in analogy to rank one spaces. In particular, the periods give rise to the translation vector of the first factor of the parallel set if $\iota=\id$.

\subsection*{Geometric interpretation of the cross ratio}

Let $x,z\in\Ftm$ and $y,w\in\Ftim$ with $x,z \op y,w$. Pick $c_x,c_z,d_y,d_w,d_w'\in\Ftsm$ such that $x$ is a face of $c_x$ and accordingly the other chambers and that $c_x\op d_y,d_w$ as well as $c_z \op d_y,d'_w$. Then we use the following notations for the horospherical retracts $\rho_{x}:=\rho_{c_x,d_y}$, $\rho_{w}:=\rho_{d_w,c_x}$, $\rho_{z}:=\rho_{c_z,d'_w}$ and $\rho_{y}:=\rho_{d_y,c_z}$. 

\bl\label{lem geometric inter cross ratios} Let $(x,y,z,w)\in\Atop$ and let $\rho_{x},\rho_{w},\rho_{z}$ and $\rho_{y}$ as above. Moreover, be $o$ in the unique affine apartment joining $c_x$ and $d_y$. Then for all $i\in J_{\tm}$ we have
$2\cro_{\xi_i}(x,y,z,w)= b_{x_{\xi_i}}(o,\rho_x \rho_w \rho_z \rho_y(o))$.
\el

\begin{proof}
Denote by $A_{xy}$ the unique affine apartment joining $c_x$ and $d_y$. Then $\rho_{d_y,c_x}$ restricted to $A_{xy}$ is the identity, i.e. $\rho_{d_y,c_x}(o)=o$. Therefore Proposition \ref{prop Gromov product as Busemann function} implies that $2(x|y)_{o,\xi_i}= b_{x_{\xi_i}}(o, o)=0$. 

By definition $\rho_y(o)$ is contained in the unique affine apartment joining $c_z$ and $d_y$. Then in the same way it follows that $(z|y)_{\rho_y(o),\xi_i}=0$. Moreover, equation \eqref{eq horospherical retract and busemann} yields $b_{y_{\iota\xi_i}}(o,\rho_y(o))=b_{y_{\iota\xi_i}}(o, o)=0$.

We can use Proposition \ref{proposition basepoint change} and again equation \eqref{eq horospherical retract and busemann} to derive that
\begin{align*}
2(z|y)_{o,\xi_i}=2(z|y)_{\rho_y(o),\xi_i}+ b_{z_{\xi_i}}(o,\rho_y(o))+ b_{y_{\iota\xi_i}}(o,\rho_y(o))= b_{z_{\xi_i}}(o,\rho_z\rho_y(o)).
\end{align*}

In a very similar way we get
\begin{align*}
2(z|w)_{o,\xi_i} & = b_{z_{\xi_i}}(o,\rho_z \rho_y(o))+ b_{w_{\iota\xi_i}}(o,\rho_w \rho_z\rho_y(o))\\
2(x|w)_{o,\xi_i} & =b_{x_{\xi_i}}(o,\rho_x\rho_w \rho_z \rho_y(o))+ b_{w_{\iota\xi_i}}(o, \rho_w \rho_z \rho_y(o)).
\end{align*}

Using that $\cro_{\xi_i}(x,y,z,w)= -(x|y)_{o,\xi_i}-(z|w)_{o,\xi_i}+(x|w)_{o,\xi_i}+(z|y)_{o,\xi_i}$, we get $2\cro_{\xi_i}(x,y,z,w)= b_{x_{\xi_i}}(o,\rho_x \rho_w \rho_z \rho_y(o))$.
\end{proof}

\bp Let $\rho_{x},\rho_{w},\rho_{z}$ and $\rho_{y}$ as before. Let $o$ be in the unique affine apartment joining $c_x, d_y$ such that we have under the identification $P(x,y)\simeq \am_{\tm}\times CS(x,y)$ that $\pi(o)=0\in\am_{\tm}$, where $\pi$ is the projection to the first factor (also assume $x\simeq \at^{+}$). Then
$2\cro_{\tm}(x,y,z,w)=\pi (\rho_x \rho_w \rho_z \rho_y(o))$.
\ep 

\begin{proof}
By construction we have that $o, \rho_x \rho_w \rho_z \rho_y(o)$ are in the unique affine apartment joining $c_x$ and $d_y$. Then by Lemma \ref{lem horospheres in parallel set} and from similar arguments as around equation \eqref{eq gromov product as busemann in apartment} we can derive that $b_{x_{\xi_i}}(o,\rho_x \rho_w \rho_z \rho_y(o))=\langle \xi_i,\pi (\rho_x \rho_w \rho_z \rho_y(o))\rangle$ for all $i\in J_{\tm}$. Together with Lemma \ref{lem geometric inter cross ratios} and the definition of $\cro_{\tm}$ we get
\begin{align*}
2\cro_{\tm}(x,y,z,w)=  \sum_{i\in J_{\tm}} \langle \xi_i,\pi (\rho_x \rho_w \rho_z \rho_y(o))\rangle \alpha^{\tm}_i.
\end{align*}
The $\xi_i\in \at$ for $i\in J_{\tm}$ form a basis of $\at$. Moreover, for all $i\in J_{\tm}$ we have that $\langle 2\cro_{\tm}(x,y,z,w), \xi_i\rangle= \langle \xi_i, \pi (\rho_x \rho_w \rho_z \rho_y(o))\rangle$. Thus it follows that $2\cro_{\tm}(x,y,z,w)=\pi (\rho_x \rho_w \rho_z \rho_y(o))$.
\end{proof}

\section{Cross ratio preserving maps}\label{sec cross ratio preserving maps}

We assume in this section that \textbf{$\tm$ is self opposite}, i.e. $\tm=\iota \tm$.

\bd
Let $M_i$, $i=1,2$ be either both symmetric spaces or thick Euclidean buildings. A map $f:\Ftmo (M_1)\to \Ftmt(M_2)$ is called \emph{$\xi_1$-Moebius map} (or \emph{cross ratio preserving}) if there exists $\xi_i\in\inn(\tm_i)$ such that 
$\cro_{\xi_1}(x,y,z,w)=\cro_{\xi_2}(f(x),f(y),f(z),f(w))$
for all $(x,y,z,w)\in \Ato$ - we in particular assume that $f(\Ato)\subset \Att$.
\ed

If $f$ is a $\xi_1$-Moebius map with respect to $\xi_1,\xi_2$, we also denote this by $\cro_{\xi_1}=f^{*}\cro_{\xi_2}$. If $\xi_1$ is clear out of the context, we sometimes call $f$ just \emph{Moebius map}. Moreover, for any map $f:\Ftmo (M_1)\to \Ftmt(M_2)$ we denote $f^{*}\cro_{\xi_2}(x,y,z,w):=\cro_{\xi_2}(f(x),f(y),f(z),f(w))$ for $x,y,z,w\in \Ftmo(M_1)$.

\bl
Let $x,y\in\Ftm$. Then there exists $z\in\Ftm$ with $z\op x,y$.
\el

\begin{proof}
We take $c_x,c_y\in \Ftsm$ such that $x$ is a face of $c_x$ and $y$ is a face $c_y$. Then there exists $c_z\in\Ftsm$ with $c_z\op c_x,c_y$ \cite[5.1]{Abramenko-vanMaldeghem}. Be $z$ the face of $c_z$ which is of type $\tm$. Then $z\in\Ftm$ with $z\op x,y$.
\end{proof}

\bl Let $f:\Ftmo(M_1)\to \Ftmt(M_2)$ be a $\xi_1$-Moebius map. Then for $x,y\in \Ftmo(M_1)$ we have that $x\op y$ if and only if $f(x)\op f(y)$.
\el  

\begin{proof}
Let $x,y\in \Ftmo(M_1)$ be given. Choose $z_1,z_2,z_3\in \Ftmo(M_1)$ such that $z_3 \op x$; $z_2\op y,z_3$ and $z_1\op x,z_2$. From Corollary \ref{corollary-not-opposite-equivalent-distance-zero} we know that $\cro_{\xi_1}(x,y,z_2,z_3)=r$ and $\cro_{\xi_1}(x,z_1,z_2,z_3)\neq \pm\infty$, i.e. $x\op y\Longleftrightarrow r\neq -\infty$. Since $\cro_{\xi_1}=f^{*}\cro_{\xi_2}$, we can derive that $f(z_2)\op f(z_3)$ and therefore we have $f(x)\op f(y)\Longleftrightarrow r\neq -\infty$. In particular, $f(x)\op f(y)\Longleftrightarrow x\op y$.
\end{proof} 

A map $f:\Ftmo(M_1)\to \Ftmt(M_2)$ such that for all $x,y\in \Ftmo(M_1)$ it holds that $x\op y$ if and only if $f(x)\op f(y)$ is called \emph{opposition preserving}.

\bl  
Let $f:\Ftmo(M_1)\to\Ftmt(M_2)$ be a $\xi_1$-Moebius map. Then $f$ is injective.
\el

\begin{proof}
Assume there exist $x\neq y\in\Ftmo(M_1)$ with $f(x)=f(y)$. Take $a\in\Ftmo(M_1) $ with $a\op x$ and $a\nop y$: For example take an apartment which contains $x$ and $y$. Take $a$ opposite of $x$ in this apartment. Then $x\neq y$ implies that $a\nop y$ - opposite points are unique in apartments.
 
In addition, choose $z,w\in\Ftmo(M_1)$ such that $z\op a$ and $w\op z,x$. Then $\cro_{\xi_1}(x,a,z,w)\neq \pm \infty$ and $\cro_{\xi_1}(y,a,z,w)=-\infty$ or is not defined; but
\begin{align*}
\cro_{\xi_1}(x,a,z,w)=f^{*}\cro_{\xi_2}(x,a,z,w)=f^{*}\cro_{\xi_2}(y,a,z,w)= \cro_{\xi_1}(y,a,z,w),
\end{align*}
contradicting $\cro_{\xi_1}(x,a,z,w)\neq \cro_{\xi_1}(y,a,z,w)$. Hence $f(x)\neq f(y)$ if $x\neq y$.
\end{proof}

\bd A surjective $\xi_1$-Moebius map is called a \emph{$\xi_1$-Moebius bijection}.
\ed

When restricting to the full flag space we can apply the following result due to Abramenko and van Maldeghem.\footnote{We remark that every spherical building is \emph{2-spherical} as in the notation of \cite{Abramenko-vanMaldeghem}. Moreover, the buildings at infinity of symmetric spaces and thick Euclidean building are thick - hence we can apply their result.}

\bp\label{pro cr preserving extends to building automorphism} (Corollary 5.2 of \cite{Abramenko-vanMaldeghem})
Let $f:\Ftsm(M_1)\to\Ftsm(M_2)$ be a surjective map that preserves opposition. Then $f$ extends in an unique way to an automorphism of the building $f:\Di M_1 \to \Di M_2 $.
\ep 

\bl \label{lem build iso of join of spher buildings}
Let $B=B_1\circ \ldots\circ B_k$ and $B'=B'_1\circ \ldots\circ B'_{k'}$ be joins of irreducible thick spherical buildings. Moreover, be $f: B\to B'$ a building isomorphism. Then $k=k'$ and there exists a permutation $\mathrm{s}$ on $k$ numbers such that $f=f_1\times\ldots\times f_k$ with $f_i:B_i\to B'_{\mathrm{s}(i)}$ building isomorphisms.
\el

\begin{proof}
That $f$ is a building isomorphism implies that $B$ and $B'$ are modeled over the same spherical Coxeter complex, i.e. over the Coxeter complex to $W=W_1\times\ldots\times W_k$, where $W_i$ are irreducible Coxeter groups. The irreducibility of the buildings $B_i,B'_i$ yields then that $k=k'$.

Assume without loss of generality that $|W_1|\leq |W_i|$ for all $i=1,\ldots,k$. Let $x_1$ be a chamber in $B_1$. Then $x_1$ is a simplex in $B$. We know that $\Res(x_1)$ is a spherical building over the spherical Coxeter complex to $W_2\times\ldots\times W_k$. As $f$ is a building isomorphism, we derive that $f(\Res(x_1))=\Res(f(x_1))$ is a spherical building over $W_2\times\ldots\times W_k$. If $f(x_1)$ would not correspond to a chamber in an irreducible factor $B'_i$, then there would be a subgroup $W'$ of $W$ isomorphic to $W_2\times\ldots\times W_k$ such that the projection of $W'$ to each $W_i$ is non-trivial (as $W_1$ is minimal). 
This would yield a decomposition of $W_2\times\ldots\times W_k$ into $k$ Coxeter groups, which contradicts the irreducibility of the factors. In particular, up to reordering $\Res(f(x_1))$ is a spherical building over $W_1\times W_3\times\ldots\times W_k$ and $W_1$ is isomorphic to $W_2$. Thus $f(x_1)=y_2$ for a chamber $y_2\in B'_2$. 
Since $f$ is a building isomorphism it maps all simplices of the same type as $x_1$ to simplices of the same type as $y_2$ i.e. it maps the chambers of $B_1$ to chambers of $B_2'$. In particular, $f$ induces a building isomorphism $f_1=f_{|B_1}:B_1\to B_2'$ ($B_1$ is naturally a subset of $B$, namely the set of simplices of $B$ fully contained in $B_1$) and thus $f=f_1\times f_0$ for a building isomorphism $f_0:B_2\circ \ldots\circ B_k\to B'_1\circ B_3' \ldots\circ B'_{k}$. A straight forward induction yields the result.
\end{proof}

We remark that multiplying the metric of a space $M$ by some positive constant $\alpha$, yields that the Gromov product on $\Ftm(\alpha M)$ is given by $(\cdot \:|\:\cdot)_{\xi,\alpha M}= \alpha (\cdot \:|\:\cdot)_{\xi,M}$ and hence also $\cro_{\xi,\alpha M}= \alpha\cro_{\xi, M}$. Moreover, there is a natural identification of $\Ftm(\alpha M)$ with $\Ftm(M)$.

\bl\label{lem moebius bij in products}
Let $M_i=M_i^1\times\ldots\times M_i^k$ be products of either irreducible symmetric spaces or irreducible thick Euclidean buildings. Moreover, be $f:\Ftsm(M_1)\to\Ftsm(M_2)$ a $\xi_1$-Moebius bijection. Then there exists a permutation $\mathrm{s}$ on $k$ numbers such that $f=f_1\times\ldots\times f_k$ with $f_i:\Ftsm(\hat{M}_1^i)\to \Ftsm(M_2^{\mathrm{s}(i)})$ a $\xi_1^{i}$-Moebius bijection and $\hat{M}_1^i$ is the space $M_1^i$ with its metric rescaled (for the types $\xi_1^{i}$ see the proof).
\el

\begin{proof}
Let $f:\Di M_1 \to \Di M_2$ be the building isomorphism from Proposition \ref{pro cr preserving extends to building automorphism}. From Lemma \ref{lem build iso of join of spher buildings} we get a permutation $\mathrm{s}$ on $k$ letters and building isomorphisms $f_i:\Di M_1^i\to \Di M_2^{\mathrm{s}(i)}$ such that
\begin{align*}
 f=f_1\times\ldots\times f_k:\Di M_1^1 \circ\ldots\circ \Di M_1^k \to \Di M_2^{\mathrm{s}(1)}\circ\ldots\circ \Di M_2^{\mathrm{s}(k)}.
\end{align*}
Moreover, we know from Proposition \ref{pro cross ratio as product} that $\cro_{\xi_i}
= \mu_i^1 \cro_{\xi_i^1}+\ldots+ \mu_i^k \cro_{\xi_i^k}$ with $\xi_i^j\in \sm_i^j$ for $i=1,2$ and $j=1,\ldots,k$ and $\mu_i\in S_k^{+}$ such that $\xi_i=\pi_i(\xi_i^1,\ldots,\xi_i^k,\mu_i)$ with $\pi_i$ as in the proposition (the numbers in the exponent are for indexing, not powers). Fix $(x_0,y_0,z_0,w_0)\in \mathrm{Flag}_{\sm_2} (M_1^2) \circ\ldots\circ \mathrm{Flag}_{\sm_k} (M_1^k)$ with $x_0,z_0 \op y_0,w_0$. Then for any $(x_1,y_1,z_1,w_1)\in \mathcal{A}_{\sm_1}$ we get
\begin{align*}
&\mu_1^1 \cro_{\xi_1^1}(x_1,y_1,z_1,w_1)+(\mu_1^2 \cro_{\xi_1^2}\ldots+ \mu_1^k \cro_{\xi_1^k})(x_0,y_0,z_0,w_0)\\
=~ &\mu_2^{\mathrm{s}(1)} f_1^{*} \cro_{\xi_2^{\mathrm{s}(1)}}(x_1,y_1,z_1,w_1)+f_0^{*}(\mu_2^{\mathrm{s}(2)} \cro_{\xi_2^{\mathrm{s}(2)}}\ldots+ \mu_2^{\mathrm{s}(k)} \cro_{\xi_2^{\mathrm{s}(2)}})(x_0,y_0,z_0,w_0)
\end{align*}
with $f_0=f_2\times\ldots\times f_k$. The equality also holds when we replace $(x_0,y_0,z_0,w_0)$ with $(z_0,y_0,x_0,w_0)$. Moreover, we have $(\mu_1^2 \cro_{\xi_1^2}\ldots+ \mu_1^k \cro_{\xi_1^k})(x_0,y_0,z_0,w_0)=-(\mu_1^2 \cro_{\xi_1^2}\ldots+ \mu_1^k \cro_{\xi_1^k})(z_0,y_0,x_0,w_0)$. Hence we derive that 
\[\mu_1^1 \cro_{\xi_1^1}(x_1,y_1,z_1,w_1) =\mu_2^{\mathrm{s}(1)} f_1^{*}\cro_{\xi_2^{\mathrm{s}(1)}}(x_1,y_1,z_1,w_1).\] 
As $(x_1,y_1,z_1,w_1)$ was arbitrary in $\mathcal{A}_{\sm_1}$ we get $\mu_1^1\cro_{\xi_1^1}=\mu_2^{\mathrm{s}(1)} f_1^{*}\cro_{\xi_2^{\mathrm{s}(1)}}$. In the same way it follows for all $i=1,\ldots,k$ that $\mu_1^i \cro_{\xi_1^i} =\mu_2^{\mathrm{s}(i)} f_i^{*}\cro_{\xi_2^{\mathrm{s}(i)}}$.

If we rescale the metric on $M_1^i$ by $\mu_2^{\mathrm{s}(i)}\slash \mu_1^i$ - denote this space by $\hat{M}_1^i$ - then $f_i:\Di \hat{M}_1^i \to \Di M_2^{\mathrm{s}(i)}$ restricts to a Moebius bijection on the chamber sets, i.e. we get a Moebius bijection $f_i:\Ftsm( \hat{M}_1^i)\to \Ftsm(M_2^{\mathrm{s}(i)})$.
\end{proof}

We will need the following fact:

\bt (\cite{Beyrer-Schroeder})\label{the BS Moebius maps in trees}
Let $T_1,T_2$ be geodesically complete trees with $|\di T_i| \geq 3$. Then every isometry from $T_1$ to $T_2$ restricted to the boundary is a Moebius bijection and every Moebius bijection $f:\di T_1\to \di T_2$ can be uniquely extended to an isometry.
\et 

Let $T$ be a rank one thick Euclidean building; in particular $T$ is a tree. Then every geodesic segment in $T$ lies in an affine apartment, i.e. in a bi-inifinite geodesic. This means that $T$ is geodesically complete (in the notation of \cite{Beyrer-Schroeder}). Moreover, by definition of thickness for rank one Euclidean buildings we have that $|\di T|\geq 3$.

We remark that $\rk(T)=1$ implies that the positive chamber of the Coxeter complex $\sm_T$ consists of a single point. Thus $\Di T= \Ftsm(T)= \di T$. Hence there is a unique Gromov product $(\cdot,\cdot)_{o_T}$ for any $o_T\in T$ on $\di T^2$ and a unique cross ratio $\cro_{T}$ on $ \mathcal{A}_T\subset \di T^4$. 

Recall that a locally compact Euclidean building with discrete translation group is called a \emph{combinatorial Euclidean building}. Moreover, given a metric realization $(B,d_B)$ of a spherical building as a $\rm{CAT}(1)$ space, the \emph{cone $E_B$ over $B$} is the quotient of $B\times[0,\infty)\slash \sim$ for the equivalence relation $(b_1,t)\sim (b_2,s)\Longleftrightarrow s=0=t$ with $b_i\in B$ and $s,t\in[0,\infty)$. The metric on $E_B$ is given by $d_{E_B}((b_1,t),(b_2,s))=s^2+t^2-2st\cos(d_B(b_1,b_2))$ 

\bp\label{pro comb Euc building moebius irred}
Let $E_1,E_2$ be irreducible thick combinatorial Euclidean buildings. Then every Moebius bijection $f:\Ftsm(E_1)\to\Ftsm(E_2)$ is the restriction of an isometry $F:\hat{E}_1\to E_2$ to the boundary where $\hat{E}_1$ is $E_1$ with its metric rescaled. If $E_1$ is not the cone over a spherical building, then $F$ is unique.
\ep

\begin{proof}
If the rank is one, then the result follows from the theorem above.

If the rank is at least 2, Struyve has shown in \cite{Struyve} that every isometry between $\di E_1$ and $\di E_2$ with respect to the Tits metric is induced by an isometry after rescaling the metric on $E_1$. The isometry is unique if $E_1$ is not the cone over a spherical building. We know that $f$ induces a building isomorphism $f:\Di E_1\to \Di E_2$ and this yields an isometry $f:\di E_1\to\di E_2$ with respect to the Tits metric when viewing simplices as subset of $\di E_i$. Hence we can apply the result of Struyve.
\end{proof}

The non-uniqueness for cones over spherical buildings arises for example as follows: Let $E_B$ be a cone over a spherical building $B$. Then the identity $\id:\Ftsm(E_B)\to \Ftsm(E_B)$ is clearly a Moebius bijection. However, every homothety of $E_B$, i.e. every map $F_\lambda:E_B\to E_B$, $(b_1,t)\mapsto(b_1,\lambda t)$ for $\lambda\in (0,\infty)$, is an isometry from $F_{\lambda}:\lambda^{2}E_B\to E_B$, where $\lambda^{2}E_B$ is the space $E_B$ with its metric rescaled by $\lambda^{2}$. In particular, every $F_{\lambda}$ extends the map $\id:\Ftsm(E_B)\to \Ftsm(E_B)$ as an isometry after rescaling the metric on the domain of $F_{\lambda}$ by $\lambda^2$.
 
\bc\label{cor comb Euc building moebius}
Let $E_1$ and $E_2$ be combinatorial Euclidean buildings and let $f:\Ftsm(E_1)\to\Ftsm(E_2)$ Moebius bijection. Then one can rescale the metric of $E_1$ on irreducible factors - denote this space by $\hat{E}_1$ - such that $f$ is the restriction of an isometry $F:\hat{E}_1\to E_2$ to the boundary. If none of the irreducible factors is a cone over a spherical building the isometry $F$ is unique.
\ec 

\begin{proof}
This follows from Lemma \ref{lem moebius bij in products} and the proposition above.
\end{proof}

\subsection{Symmetric spaces}
We want to show that the above proposition and corollary hold in a similar way for symmetric spaces. We will see that we essentially only need to show that Moebius bijections are homeomorphisms. Therefore we analyze some topological properties of Moebius bijections for the case of symmetric spaces.

In this section we only consider symmetric spaces $X$. For $r\in\R$, $\xi\in\inn(\tm)$ and $x_2,y_1,y_2\in \Ftm(X)$ we define 
\begin{align*}
B^{+}_{r,\xi}(y_1,x_2,y_2)&:=\lbrace x_1\in\Ftm(X) \; |\; (x_1,y_1,x_2,y_2)\in \Af,\; \cro_{\xi} (x_1,y_1,x_2,y_2) > r\rbrace,\nonumber \\
B^{-}_{r,\xi}(y_1,x_2,y_2)&:=\lbrace x_1\in\Ftm(X) \; |\; (x_1,y_1,x_2,y_2)\in \Af,\; \cro_{\xi} (x_1,y_1,x_2,y_2) < r\rbrace.
\end{align*}
Those sets are open by the continuity of $\cro_{\xi}$ and the fact that $\Af$ is open. However, it can happen that they are empty - which holds if $x_2\nop y_1,y_2$.

\bp\label{prop-cross ratio generating topology}
Let $X$ be a symmetric space. The sets $B^{-}_{r,\xi}(y_1,x_2,y_2)$ varying over all $r\in \R$ and all $x_2,y_1,y_2\in \Ftm$ form a subbase of the topology on $\Ftm(X)$
\ep 

\begin{proof}
As mentioned, those sets are open. Thus it is enough to show that any open neighborhood $U$ of a point $x\in\Ftm(X)$ contains an open neighborhood $V$ which can be written as a finite intersection of sets of the form $B^{-}_{r,\xi}(y_1,x_2,y_2)$.

Let $x\in\Ftm(X)$ and let any neighborhood $U$ of $x$ be given. We set $K:=\Ftm\backslash U$. Then $K$ is compact and $x\notin K$.

For any $a\in K$, choose $y_a\in\Ftm(X)$ such that $y_a\op a$ and $y_a\nop x$. In addition, choose $w_a,z_a\in\Ftm(X)$ such that $w_a\op a,x$ and $z_a\op y_a,w_a$. This yields $\cro_{\xi} (x,y_a,z_a,w_a)=-\infty$ and $\cro_{\xi} (a,y_a,z_a,w_a)>r_a$ for some $r_a\in\R$ and hence $x\in B^{-}_{r_a,\xi}(y_a,z_a,w_a)$, $x\notin B^{+}_{r_a,\xi}(y_a,z_a,w_a)$, $a\in B^{+}_{r_a,\xi}(y_a,z_a,w_a)$.

Varying over all $a\in K$ the sets $B^{+}_{r_a,\xi}(y_a,z_a,w_a)$ cover $K$ and by compactness we find a finite number of points $a_i\in K$, $i=1,\ldots,l$ such that the according sets already cover $K$. We set $V:=\bigcap_{a_i:i=1,\ldots,l} B^{-}_{r_{a_i},\xi}(y_{a_i},z_{a_i},w_{a_i})$.
As a finite intersection of open sets, $V$ is open. Furthermore, $x\in V$ and hence $V$ is non-empty. By construction $V\subset K^C$ and hence $V\subset U$.
\end{proof}

\bl \label{lem moebius bijection is a homeo}
Let $f:\Ftmo(X_1)\to\Ftmt(X_2)$ be a $\xi_1$-Moebius bijection. Then $f$ is a homeomorphism.
\el

\begin{proof}
Since $f$ leaves the cross ratio invariant and is a bijection, it is immediate that $f(B^{-}_{r,\xi_1}(y,z,w))=B^{-}_{r,\xi_2}(f(y),f(z),f(w))$. This means that $f$ yields a bijection of subbases of the topology and hence $f$ is a homeomorphism.
\end{proof}

As mentioned, for a symmetric space $X$ the boundary $\Ftm(X)$ can be identified homeomorphically with $G\slash P_x$ for $P_x=\st(x)$ and $x\in \Ftm(X)$. Hence $\Ftm(X)$ can be given the structure of compact connected manifold (without boundary) - inherited from $G\slash P_x$. Using this there is a different way to characterize Moebius bijections captured in the following lemma.

\bl 
Let $X_1, X_2$ be symmetric spaces such that $\dim \Ftmo(X_1) = \dim \Ftmt(X_2)$ and $f:\Ftmo(X_1)\to\Ftmt(x_2)$ be a continuous $\xi_1$-Moebius map. Then $f$ is a homeomorphism, in particular $f$ is a $\xi_1$-Moebius bijection.
\el

\begin{proof}
Since $f$ is a $\xi_1$-Moebius map and hence injective, we know that $f:\Ftmo(X_1)\to \mathrm{Im}(f)$ is a bijection, with $\mathrm{Im}(f)$ denoting the image. Moreover, $f^{*}\cro_{\xi_2}=\cro_{\xi_1}$ implies $f(B^{-}_{r,\xi_1}(y,z,w)) =B^{-}_{r,\xi_2}(f(y),f(z),f(w))\cap \Ima(f)$.
Then Proposition \ref{prop-cross ratio generating topology} yields that $f$ maps a subbase of the topology on $\Ftmo(X_1)$ into a subbase of the topology on $\mathrm{Im}(f)$ equipped with the subset topology. Hence $f:\Ftmo(X_1)\to \mathrm{Im}(f)$ is open and therefore a homeomorphism.

We derive that $\Ima(f)$ is compact connected submanifold of $\Ftmt(X_2)$ of the same dimension. However, $\Ftmt(X_2)$ is a compact connected manifold without boundary and hence the only such submanifold is $\Ftmt(X_2)$ itself, i.e. $\Ima(f)=\Ftmt(X_2)$ - which proves the claim.
\end{proof}

\bt\label{thm Moebius maps in sym space}
Let $X_1,X_2$ be symmetric spaces of rank at least two with no rank one de Rham factors and let $f:\Ftsm(X_1)\to\Ftsm(X_2)$ be a $\xi_1$-Moebius bijection. Then one can multiply the metric of $X_1$ by positive constants on de Rham factors - denote this space by $\hat{X}_1$ - such that $f$ is the restriction of an unique isometry $F:\hat{X}_1\to X_2$ to $\Ftsm(X_1)$. 
\et

\begin{proof}
We know that a $\xi_1$-Moebius bijection $f:\Ftsm(X_1)\to\Ftsm(X_2)$ can uniquely be extended to a building isomorphism $f:\Di X_1\to \Di X_2$. Moreover, $f$ is a homeomorphism on the chamber sets $\Ftsm(X_i)$ by Lemma \ref{lem moebius bijection is a homeo}. Then for such maps the result is known \cite[Sc.3.9]{Eberlein}.
\end{proof}

Actually all we need for the above result is that $f:\Ftsm(X_1)\to\Ftsm(X_2)$ is opposition preserving and a homeomorphism. However, when dealing also with rank one factors we really need Moebius maps. 

\bc\label{cor moebius bij for sym spaces}
Let $X_1$ and $X_2$ be symmetric spaces of non-compact type and let $f:\Ftsm(X_1)\to\Ftsm(X_2)$ be a Moebius bijection. Then one can rescale the metric of $X_1$ on de Rham factors - denote this space by $\hat{X}_1$ - such that $f$ is the restriction of an unique isometry $F:\hat{X}_1\to X_2$ to the boundary.
\ec 

\begin{proof}
This follows from Lemma \ref{lem moebius bij in products} together with the theorem above and the fact that Moebius bijections of rank one symmetric spaces can be uniquely extended to isometries. For the latter result see \cite{Bourdon}.
\end{proof}

\subsection{Rescaling on irreducible factors}

In this generality it is not possible to drop the scaling on the irreducible factors in the Corollaries \ref{cor comb Euc building moebius}, \ref{cor moebius bij for sym spaces} and Theorem \ref{thm Moebius maps in sym space}. For example consider the following situation:
Let $M_0$ be a symmetric space or a combinatorial Euclidean building. We set $M_1:=\mu^{-1}_1 M_0$, $M_2:=\mu^{-1}_2 M_0$ for $\mu_i>0$ with $\mu_1^2+\mu_2^2=1$ and $M:=M_1\times M_2$ - here $M_i=\mu^{-1}_i M_0$ means we take the space $M_0$ with its metric multiplied by $\mu_i^{-1}$. Moreover, we define $f:\Ftsm(M)\to \Ftsm(M)$ by $f(x,y):=(y,x)$.

Let $\xi\in \inn(\sm_0)$ and $\sm_0$ the fundamental of the space $M_0$. Consider the cross ratio $\cro_{\pi(\xi,\xi,(\mu_1,\mu_2)),M}=\mu_1\cro_{\xi,M_1}+\mu_2\cro_{\xi,M_2}$ - cp. Proposition \ref{pro cross ratio as product}. As mentioned, we have $\mu_1\cro_{\xi,M_1}=\cro_{\xi,M_0}=\mu_2\cro_{\xi,M_2}$ and hence $f$ is a $\pi(\xi,\xi,(\mu_1,\mu_2))$-Moebius bijection.

We see that $f$ is induced by a map $F:=F_1\times F_2:M_1\times M_2\to M_2\times M_1$, such that $F_i:\Ftsm(M_i)\to \Ftsm(M_j)$, $i\neq j$ is the identity (under the natural identification with $\Ftsm(M_0)$). As $F$ and hence the $F_i$ shall be isometries, it follows that $F(p,q)=(q,p)$ and clearly $F$ is an isometry only after rescaling on de Rham factors.

Let $M_1$ be a symmetric space or a combinatorial Euclidean building and assume that the image of $\cro_{\sm,M_1}$ lies not in a proper subspace of $\am_{M_1}$. Then the above situation is essentially the only possibility where rescaling can appear:

Let $M_1,M_2$ be irreducible. In addition, be $f:\Ftsm(M_1)\to\Ftsm(M_2)$ a $\xi_1$-Moebius bijection, i.e. $\cro_{\xi_1}=f^{*}\cro_{\xi_2}$. Then we know that we can rescale the metric on $M_1$ by some positive number $\mu_1$, such that $f$ is induced by an isometry $F:\mu_1 M_1\to M_2$. Thus Proposition implies \ref{prop isometries are Moebius maps} $f^{*}\cro_{\xi_2}=\cro_{\xi'_1,\mu_1 M_1}=\mu_1 \cro_{\xi'_1,M_1}$ for $\xi'_1\in\sm_1$ with $F_{\sm}(\xi'_1)=\xi_2$.

However, it follows from the assumption on $\cro_{\sm,M_1}$ together with Lemma \ref{lem vector valued cro gives usual by inner product} that $\cro_{\xi}\neq \alpha\cro_{\xi'}$ for $\xi\neq \xi'\in\sm_1$ and any $\alpha\in\R$. Therefore $\cro_{\xi_1,M_1}=f^{*}\cro_{\xi_2}=\mu_1 \cro_{\xi'_1,M_1}$ implies $\xi_1=\xi'_1$ and $\mu_1=1$ - in particular $f$ is induced by an isometry without rescaling the metric.

We remark that for symmetric spaces with $\iota=\id$ the image of $\cro_{\sm}$ is all of $\am$. This follows from the fact that every vector of $\am$ can be realized as a translation vector of a hyperbolic element in $G$. Then the periods of those elements in $G$ are exactly those translation vectors, as seen in Proposition \ref{pro periods of isometries}. Hence the above discussion applies.

\bc
Let $M$ either be a symmetric space or a combinatorial Euclidean building with none of the irreducible factors being a cone over a spherical building. In addition, assume that the image of $\cro_{\sm}$ is not contained in a proper subspace of $\am$. Let $\xi_0\in\sm$ be the center of gravity of $\sm$. Then there is a one-to-one correspondence between $Iso(M)$ and $\xi_0$-Moebius bijections.
\ec

\begin{proof}
Let $g\in Iso(M)$ and $g_{\sm}:\sm\to\sm$ the induced map. Then $g_{\sm}$ is an isometry with respect to the angular metric, hence $g_{\sm}$ stabilizes the center of gravity $\xi_0$ of $\sm$. Therefore Proposition \ref{prop isometries are Moebius maps} yields a $\xi_0$-Moebius bijections for each $g\in Iso(M)$.

On the other hand, by Corollaries \ref{cor comb Euc building moebius} and \ref{cor moebius bij for sym spaces}, we know that each $\xi_0$-Moebius bijections is induced by a unique isometry - after possible rescaling on irreducible factors. However, following the above discussion we can exclude rescaling of the metric: 

Let $f$ be a $\xi_0$-Moebius bijection and let $f=f_1\times\ldots\times f_k$ be the decomposition on irreducible factors $M_1,\ldots, M_k$ as in Lemma \ref{lem moebius bij in products}. 
Assume w.l.o.g. that $f_1:\Ftsm(M_1)\to \Ftsm(M_2)$, i.e. $M_1, M_2$ are isometric after possibly rescaling the metric.
From Proposition \ref{pro cross ratio as product} we know $\cro_{\xi_0}=\mu_1 \cro_{\xi_1,M_1}+\mu_2 \cro_{\xi_2,M_2}+\ldots+\mu_k \cro_{\xi_k,M_k}$.
However, $\xi_0\in\sm$ being the center of gravity of $\sm$ and $M_1, M_2$ isometric after possibly rescaling the metric implies $\mu_1=\mu_2$ and $\xi_1\simeq \xi_2$. Then $f_1$ is $\xi_1$-Moebius bijection between irreducible spaces. From the above discussion it follows that it is induced by an isometry without rescaling the metrics. The same argument implies the result for all $f_i$ and hence the claim follows.
\end{proof}

\subsection{General Euclidean buildings}

In this section we consider general Euclidean buildings, i.e. in particular non-locally compact ones. The goal is again to show that Moebius bijections are induced by isometries. However, now we will need the vector valued cross ratio $\cro_{\sm}$ to derive such a result.

Let $E$ be a thick Euclidean building considered with the complete apartment system. Let $x\in \Ftm(E)$ and $y\in\Ftim(E)$ with $x\op y$ and $\tm$ is a codimension 1 face of $\sm$ - in this case $x,y$ are called \emph{panels} of the building $\Di E$. Then metrically we have the splitting $P(x,y)\simeq\at\times CS(x,y)$, where $ CS(x,y)$ is a Euclidean building of rank $\mathrm{rk}(E)-\dim \at =1$, i.e. $CS(x,y)$ is an $\R$-tree. This tree is called \emph{wall tree} and will be denoted by $T_{xy}$. One can show that the isomorphism type of $T_{xy}$ does not depend on the choice of $y\in \Ftim(E)$ with $y\op x$ \cite{Kramer-Weiss}; hence the isomorphism class of $T_{xy}$ will be denoted by $T_x$.

We recall that the \emph{residue} of an element $z\in \Di E$ is defined by $\Res(z)=\{ w\in \Di E~ | z\subsetneq w\}$. In case of a panel $x\in \Di E$ we have that  $\Res(x)$ consists of all the chambers in $\Di E$ containing $x$. 

It is known that one can naturally identify $\Res(x)\simeq\di T_x$. Let us describe this identification: 
Fix $y\op x$ and consider $T_{xy}$ in the isomorphism class $T_x$. Let $o\in P(x,y)$. Then we can identify $P(x,y)\simeq\at\times T_{xy}$ such that $o\simeq(0,o_T)$ and $x\simeq \di \at^{+}$ - recall that $\at^{+}=\at\cap \am^{+}$. Then there is a one-to-one correspondence of chambers in $\Res(x)$ with (specific) Weyl sectors in $P(x,y)$ with tip $o$ \cite[Cor.~1.9.]{Parreau}.\footnote{Here \emph{Weyl sector} includes also all translates in an affine apartment of the Weyl sectors we have considered so far.}
The affine apartments in $P(x,y)\simeq\at\times T_{xy}$ containing $o$ are of the form $\at\times \gamma$, where $\gamma$ is a bi-infinite geodesic ray in $T_{xy}$ passing through $o_T$ (those are easily seen to be isometric to $\R^r$). 
By definition every Weyl sector is contained in an affine apartment; hence we can derive that every Weyl sector with tip $o$ and boundary chamber $c\in\Res(x)$ is contained in $\at^{+}\times \gamma_{o_Tz}$ where $\gamma_{o_Tz}$ is a geodesic ray in $T_{xy}$ from $o_T$ to a boundary point $z\in\di T_{xy}$.
This yields a one-to-one correspondence of $\Res(x)$ with geodesic rays emanating from $o_T$. As those rays are in one-to-one correspondence with $\di T_{xy}$; hence we get $\Res(x)\simeq\di T_x$ as claimed.

\begin{Remark}\label{rem residues and wall tree identification}
It follows that for $z\in \di T_{xy}, c\in\Res(x)$ and $d\in \Res(y)$ we have that $z\simeq c$ and $z\simeq d$ under $\Res(x)\simeq\di T_{xy}$, $\Res(y)\simeq\di T_{xy}$ respectively if and only if the Weyl sectors with tip $o=(0,o_T)$ defining $c,d$ are contained in $\at^{+}\times \gamma_{o_Tz}$, $\at^{-}\times \gamma_{o_Tz}$, respectively.
\end{Remark}

By definition $\Res(x)$ is the set of chambers that contain $x$. Hence there is a unique corner $\xi_{x}$ of $\sm$ such that $c_{\xi_x}\notin x$ for every chamber $c\in\Res(x)$. In the same way we get a type from $y$ and it is immediate that this type equals $\iota\xi_x$ - following for example from the fact that $x\in\Ftm$ implies that $y\in \Ftim$.

\bl\label{lem gromov product and wall trees}
Let $x, y$ be opposite panels in $\Di E$ and $T_{xy}$ the associated tree. Let $z_c,z_d\in\di T_{xy}$, $c\in \Res(x)$ such that $c\simeq z_c$ under $\Res(x)\simeq\di T_{xy}$ and $d\in \Res(y)$ such that $d\simeq z_d$ under $\Res(y)\simeq\di T_{xy}$. Then $(c|d)_{o,\xi_x}=\sin(\alpha)(z_c|z_d)_{o_T}$ where $o\simeq (0,o_T)$ under $P(x,y)\simeq\at\times T_{xy}$ and $\alpha\in (0,\pi)$ does only depend on $\sm$ and the type of $x$.
\el

\begin{proof}
Let $\gamma_c,\gamma_d$ be the geodesics in $P(x,y)$ from $o$ to $c_{\xi_{x}}$ and $d_{\xi_{y}}$, respectively. The splitting $P(x,y)\simeq\at\times T_{xy}$ yields geodesics $\gamma_{x}, \gamma_{y}$ in $\at$  from $0$ and $\gamma_{z_c}, \gamma_{z_d}$ in $T_{xy}$ eminating from $o_T$ such that $\gamma_c(t)=(\gamma_x(t),\gamma_{z_c}(t))$ and $\gamma_d(t)=(\gamma_y(t),\gamma_{z_d}(t))$ - while $\gamma_c, \gamma_d$ are unit speed, the geodesics $\gamma_{x}, \gamma_{y}, \gamma_{z_c}$ and $\gamma_{z_d}$ are not. It is clear that the geodesics $\gamma_{x}, \gamma_{y}$ do not depend on the choice of $c,d$ and are in opposite directions (since the $\gamma_c, \gamma_d$ are): 
The geodesics $\gamma_c,\gamma_d$ are along those corners of Weyl sectors that are not contained in $\at$. Since Weyl sectors are isometric to convex subsets of $\R^r$, it reduces to Euclidean geometry; for example $\gamma_{x}$ is the geodesic in $\at$ from $0$ to the point in $x$ of type $\pi_{\tau_x}(\xi_x)$, where $\pi_{\tau_x}$ is the orthogonal projection from $\sm$ to $\tau_x$ and $\tau_x$ is the type of $x$. 

Let from now on $\gamma_{x}, \gamma_{y}, \gamma_{z_c}$ and $\gamma_{z_d}$ be the geodesics as above but now parametrized such that they are unit speed. Then the above discussion implies that $d(\gamma_{x}(t),\gamma_{y}(t))=2t$. 
Let $\alpha$ be the angle of $\xi_x$ and $\pi_{\tau_x}(\xi_x)$. Then we have $\gamma_c(t)=(\gamma_x(\cos(\alpha) t),\gamma_{z_c}(\sin(\alpha)t))$. Basic facts of trees imply that $d(\gamma_{z_c}(t),\gamma_{z_d}(t))=2t-2(z_c|z_d)_{o_T}$ for $t\geq (z_c|z_d)_{o_T}$ - see e.g. \cite{Beyrer-Schroeder}. Altogether,
\begin{align*}
(c|d)_{o,\xi_x}=& \lim_{t\to\infty} t-\frac{1}{2} \sqrt{4\cos^2(\alpha) t^2+(2\sin(\alpha)t-2(z_c|z_d)_{o_T})^2}\\
=& \lim_{t\to\infty} t -\sqrt{t^2-2t\sin(\alpha)(z_c|z_d)_{o_T}+ (z_c|z_d)_{o_T}^2}=\sin(\alpha)(z_c|z_d)_{o_T},
\end{align*}
while the last equality follows from a Taylor series in the same way as we have seen several times before.
\end{proof}

\bc\label{cor cross ratio of tree as usual one}
The cross ratio on $\di T_{xy}$ is given by $\cro_{T_{xy}}(z_1,w_1,z_2,w_2)=\sin(\alpha)\cro_{\xi_x}(c_1,d_1,c_2,d_2)$ where $\xi_x\in\sm$ is the corner not contained in $\tm_x$, the type of $x$, $\alpha$ is the angle between $\xi_x$ and $\tm_x$, $c_i\simeq z_i$ under $\Res(x)\simeq\di T_{xy}$ and $d_i\simeq w_i$ under $\Res(y)\simeq\di T_{xy}$.
\ec

The thickness of $E$ means that $\Di E$ is thick and therefore for every panel $x$ we have that $|\di T_x|\geq 3$ (as $\Res(x)\simeq\di T_x$), i.e. $T_x$ is thick and geodesically complete. Therefore Theorem \ref{the BS Moebius maps in trees} implies that the whole isometry class $T_x$ has a natural cross ratio $\cro_{T_x}$.

\bd
Let $E_1, E_2$ be thick irreducible Euclidean buildings. A building isomorphism $\phi:\Di E_1\to \Di E_2$ is called \emph{tree-preserving} or \emph{ecological}, if for every panel $x\in \Di E_1$ we have that $\phi_{|\Res(x)}:\Res(x)\to \Res(\phi(x))$ is induced by an isometry $\phi_x:T_x\to T_{\phi(x)}$ - i.e. $(\phi_x)_{|\di T_x}\simeq \phi_{|\Res(x)}$ under the identification $\Res(x)\simeq\di T_x$ 
\ed

\bt (Tits,~\cite[Thm~2]{Tits}) 
Let $E_1, E_2$ be thick irreducible Euclidean buildings and $\phi:\Di E_1\to \Di E_2$ an ecological isomorphism. Then $\phi$ extends to an isomorphism, i.e. an isometry after possibly rescaling the metric on $E_1$.
\et 

In a similar way as before, we call a surjective map $f:\mathrm{Flag}_{\sm_1} (E_1)\to \mathrm{Flag}_{\sm_2}(E_2)$ such that $\cro_{\sm_1}(x,y,z,w)=f^{*}\cro_{\sm_2}(x,y,z,w)$ for all $(x,y,z,w)\in \mathcal{A}_{\sm_1}$ a \emph{$\sm_1$-Moebius bijection}. We remark that to identify the image of $\cro_{\sm_1}$ with the one of $\cro_{\sm_2}$ it is already necessary that $E_1$ and $E_2$ are modeled over the same Coxeter complex, i.e. $\sm_1\simeq\sm_2=:\sm$.

It is immediate that such a map is a $\xi_0$-Moebius map, for $\xi_0$ the center of gravity of $\sm$. We assumed $f$ to be surjective, hence $f$ is a $\xi_0$-Moebius bijection and therefore $f$ can be extended uniquely to a building isomorphism $f:\Di E_1\to \Di E_2$ by Proposition \ref{pro cr preserving extends to building automorphism}.\medskip

We recall that the affine Weyl group $\hat{W}=W\ltimes T_W$ of the Coxeter complex over which a Euclidean building is defined gives a collection of hyperplanes, namely the hyperplanes of the finite reflection group $W$ together with all its translates under $T_W$. Each hyperplane defines two half spaces which we call \emph{affine half apartments}. The image of an affine half apartment under a chart map is again called \emph{affine half apartment}. 

In spherical buildings the hyperplanes associated to the spherical Coxeter group define walls in apartments and those walls separate the apartments in two halfs, called \emph{half apartments}.
One can show that the boundary of an affine half apartment $H\subset E$ defines a half apartment in $H_{\infty}\subset \Di E$ and to every half apartment in $H_{\infty}\subset\Di E$ we find an affine half apartment $H\subset E$ which has $H_{\infty}$ as its boundary.

Now, let $f:\Di E_1\to \Di E_2$ be a building isomorphism and let $x,y$ be opposite panels. The identifications $\di T_{xy}\simeq \Res(x), \di T_{xy}\simeq \Res(y)$ together with $f_{|\Res(x)}:\Res(x)\to\Res(f(x))$, $f_{|\Res(y)}:\Res(y)\to\Res(f(y))$ induce two maps $f_x,f_y:\di T_{xy}\to \di T_{f(x)f(y)}$.

\bl \label{lem wall tree maps independent of choice of residue}
Notations as above, in particular $x,y$ are opposite panels and $f_x,f_y:\di T_{xy}\to \di T_{f(x)f(y)}$ are induced by $f_{|\Res(x)}:\Res(x)\to\Res(f(x))$, $f_{|\Res(y)}:\Res(y)\to\Res(f(y))$. Then $f_x=f_y$.
\el

\begin{proof}
Let $z\in \di T_{xy}$, i.e. $z$ is an equivalence class of geodesic rays. Every ray $\gamma_z$ in the class starting at a branching point defines an affine half apartment $\at\times \gamma_z$ in $E_1$ and thus (the equivalence class of rays) defines a half apartment $H_{\infty}\subset\Di E_1$. 
Then it follows form Remark \ref{rem residues and wall tree identification} that $c\simeq z$ with $c\in\Res(x)$ if and only if $c$ is contained in the half apartment $H_{\infty}$ and in the same way $d\simeq z$ with $d\in\Res(y)$ if and only if $d$ is contained in the half apartment $H_{\infty}$. 
By assumption, $f$ is a building isomorphism, i.e. $f(H_{\infty})\subset \Di E_2$ is a half apartment with $f(x),f(y)\in f(H_{\infty})$. The metric splitting $P(f(x),f(y))\simeq \at\times T_{f(x)f(y)}$ yields that we find an affine half apartment $\at\times \gamma_w$ with $\gamma_w$ a geodesic ray in $T_{f(x)f(y)}$ and boundary point $w\in \di T_{f(x)f(y)}$ such that the boundary of this affine half apartment is exactly $f(H_{\infty})$. By definition $f(c),f(d)\in f(H_{\infty})$. Hence from Remark \ref{rem residues and wall tree identification} we get that $f(c)\simeq w\simeq f(d)$. Therefore $f_x(z)=w$ and $f_y(z)=w$.
\end{proof}

\bt\label{thm moebius maps are ecological}
Let $E_1, E_2$ be thick irreducible Euclidean buildings. Let $f:\Ftsm(E_1)\to \Ftsm(E_2)$ be a $\sm$-Moebius bijection. Then the induced isomorphism $f:\Di E_1\to \Di E_2$ is ecological and hence can be extended to an isomorphism $F:E_1\to E_2$, i.e. an isometry after possibly rescaling the metric on $E_1$
\et

\begin{proof}
What we need to show is, given a panel $x\in \Di E_1$, the induced map $f_x:\di T_x\to \di T_{f(x)}$ is the restriction of an isometry. This implies that $f$ is ecological and therefore by the Theorem of Tits induced by an isomorphism.

We fix $y\op x$ to get a tree $T_{xy}$ in the class of $T_x$. Since we are considering isometry classes of trees, it is enough to show that $f_{xy}:\di T_{xy}\to \di T_{f(x)f(y)}$ is induced by an isometry. 

Corollary \ref{cor cross ratio of tree as usual one} implies that for $z_1,w_1,z_2,w_2\in\di T_{xy}$ and $c_1,c_2\in\Res(x)$, $d_1,d_2\in\Res(y)$ with $z_i\simeq c_i$, $w_i\simeq d_i$ there is some $\alpha\in (0,\pi)$ with
\begin{align*}
\cro_{T_{xy}}(z_1,w_1,z_2,w_2)=\sin(\alpha)\cro_{\xi_x}(c_1,d_1,c_2,d_2)=\sin(\alpha)f^{*}\cro_{\xi_x}(c_1,d_1,c_2,d_2),
\end{align*}
while the last equality follows from $f$ being a $\sm$-Moebius bijection. By construction $f_{xy}:\di T_{xy}\simeq \Res(x)\to \di T_{f(x)f(y)}\simeq \Res(f(x))$ is defined in the way that $f(c_1)\simeq f_{xy}(z_1)$ under $\di T_{f(x)f(y)}\simeq \Res(f(x))$ and similar for $c_2$. In light of Lemma \ref{lem wall tree maps independent of choice of residue} we have that $f(d_i)\simeq f_{xy}(w_i)$.
Applying again Corollary \ref{cor cross ratio of tree as usual one} this yields that $\sin(\alpha)f^{*}\cro_{\xi_x}(c_1,d_1,c_2,d_2)=f_{xy}^{*}\cro_{T_{f(x)f(y)}}(z_1,w_1,z_2,w_2)$ - we remark that the $\alpha$ is the same as before as the simplices $\sm_1$ and $\sm_2$ coincide. 
Hence $f_{xy}$ is a Moebius bijection. Since $T_{xy}$ is a geodesically complete tree and the thickness of $E_1$ implies that $|\di T_{xy}|\geq 3$ we can apply Theorem \ref{the BS Moebius maps in trees} to derive that $f_{xy}$ is induced by an isometry.
\end{proof}

\bc\label{cor Moebius bij on red Euclidean build}
Let $E_1, E_2$ be thick Euclidean buildings and moreover let $f:\Ftsm(E_1)\to\Ftsm(E_2)$ be a $\sm$-Moebius bijection. Then we can rescale the metric on the irreducible factors of $E_1$ - denote this space by $\hat{E}_1$ - such that $f$ is the restriction of an isometry $F:\hat{E}_1\to E_2$ to the boundary.
\ec

\begin{proof}
Since $f$ can be extended to a building isomorphism (as we have seen before), $f$ is opposition preserving for each type of simplex. This, together with Lemma \ref{lem vector valued cro gives usual by inner product} and $f$ being a $\sm$-Moebius bijection, yield that $f^{*}\cro_{\xi}=\cro_{\xi}$ for every type $\xi\in\sm$.

Let $\sm=\sm_1\circ \ldots\circ\sm_k$ be the decomposition of $\sm$ corresponding to the decomposition of $E_i$ into irreducible factors - the decompositions coincide as both buildings are thick and modeled over the same spherical Coxeter complex. Moreover, be $f=f_1\times\ldots\times f_k$ the decomposition from Lemma \ref{lem moebius bij in products}. 

Then $f^{*}\cro_{\xi}=\cro_{\xi}$ for all $\xi\in\sm$ implies that each $f_i$ is a $\sm_i$-Moebius bijection. Thus the above theorem yields the claim.
\end{proof}

\bc
Let $E_1, E_2$ be thick irreducible Euclidean buildings. In addition, assume that there exists a wall tree $T_x$ for a panel $x\in\Di E_1$ which has more than one branching point. Let $f:\Ftsm(E_1)\to \Ftsm(E_2)$ be a $\sm$-Moebius bijection. Then $f$ can be extended to an isometry $F:E_1\to E_2$ (without rescaling the metric).

Moreover, if $E_1$ is not a Euclidean cone over a spherical building then every wall tree has more than one branching point.
\ec

\begin{proof}
From Theorem \ref{thm moebius maps are ecological} we know that we can rescale the metric by some $\mu\in\R_{+}$ such that $f$ is induced by an isometry $F:\mu E_1\to  E_2$, where $\mu E_1$ is $E_1$ with the metric rescaled by $\mu$. Let $x\in\Di E_1$ be a panel such that the wall tree $T_x$ has more than one branching point. Then clearly the wall tree of $x\in \Di \mu E_1$ is $\mu T_x$. Let $f_x:\di T_x\to \di T_{f(x)}$ be the induced map from $f$ on the wall tree. Since $F$ restricted to the boundary is $f$, the map induced from $F$ on $\di \mu T_x$ equals $f_x$. Therefore we have $\cro_{T_x}= f_x^{*} \cro_{T_{f(x)}}=\cro_{\mu T_x}=\mu \cro_{T_x}$ (the first equality follows from $f$ being a $\sm$-Moebius bijection, the second from $f_x=F_{|\di \mu T_x}$). 

By assumption $T_x$ has two branching points. The distance of those two points can be given in terms of the cross ratio - i.e. let $p,q\in T_x$ be the branching points, then there exist $z_1,z_2,w_1,w_2\in \di T_x$ such that $d(p,q)=\cro_{T_x}(z_1,w_1,z_2,w_2)$ \cite[Lem 4.2]{Beyrer-Schroeder}. Since this distance $d(p,q)$ is non-zero, we derive from $\cro_{T_x}(z_1,w_1,z_2,w_2)=\mu \cro_{T_x}(z_1,w_1,z_2,w_2)$ that $\mu=1$. Hence $F$ is an isometry without rescaling the metric on $E_1$.

The second claim is a direct consequence of Proposition 4.21. and 4.26 in \cite{Kramer-Weiss}.
\end{proof}

The second claim of Theorem \ref{thm B} follows now from the fact that every $\sm$-Moebius bijection splits as a product of $\sm_i$-Moebius bijections on irreducible factors, as in the proof of Corollary \ref{cor Moebius bij on red Euclidean build}. The corollary above implies that those $\sm_i$-Moebius bijections induce isometries without the need of rescaling.

\section{Appendix}
Here, we determine the cross ratios that we construct explicitly for the symmetric spaces $X(n):=\text{SL}(n,\R)\slash \text{SO}(n,\R)$. We will use the notation as in Example \ref{ex Sl n}. 

The map $g\cdot \text{SO}(n,\R)\mapsto gg^t$ yields an identification of $X(n)$ with the space
$
P_n=\lbrace A\in \text{Mat}(n\times n,\mathbb{R})|A=A^t\wedge \det(A)=1\wedge A \text{ is positive definite}\rbrace$.
The action of $g\in \text{SL}(n,\R)$ on $A\in P_n$ is given by $g\cdot A=gAg^t$. 
By definition of the cross ratio, it will be enough to determine $(\cdot\,|\,\cdot)_{I_n,\lambda}$ with $I_n$ being the identity matrix in $P_n$ and $\lambda=(\lambda_1,\ldots,\lambda_l)$ be identified with some type.

Let $\tm=(i_1,\ldots,i_l)$, $i_j\in \lbrace 1,\ldots,n\rbrace$ such that $i_l=n$, $i_j<i_{m}$ for $1\leq j<m\leq l\leq n$ and be $\textsc{S}_{\tm}$ be the corresponding standard flag, i.e. $\textsc{S}_{\tm}=(V_{i_1},\ldots, V_{i_l})$ for $V_{i_j}=\text{span} \lbrace e_1,e_2,\ldots, e_{i_j}\rbrace$. 
Let $\textsc{S}_{\iota\tm}$ be the standard opposite flag to $\textsc{S}_{\tm}$, i.e. $\textsc{S}_{\iota\tm}=(V^{*}_{i_{l-1}},\ldots, V^{*}_{i_1},\R^n)$ with $V^{*}_{i_j}=\text{span} \lbrace e_n,e_{n-1},\ldots, e_{i_j+1}\rbrace$. 
Furthermore, be $\lambda=(\lambda_1,\ldots,\lambda_l)\in \mathbb{R}^l$ such that $\lambda_j>\lambda_{j+1}$, $\sum_{j=1}^l m_j \lambda_j=0$ for $m_j=\dim V_{i_j} - \dim V_{i_j -1}$ if $j>1$, $m_1=\dim V_{i_1}$ and $\sum_{j=1}^l m_j \lambda^2_j=1$.

\begin{claim}
Notations as before, $k,h\in \mathrm{SO}(n,\R)$ and denote by $\hat{h}_i$ the $i$-th column of the matrix $h$ and accordingly $\hat{k}_i$. Then
\begin{align*}
(k\textsc{S}_{\tm}|h\textsc{S}_{\iota\tm})_{I_n,\lambda}=n \sum_{j=1}^{l-1} (\lambda_{j+1}-\lambda_{j}) \log|\det (\,\hat{k}_1\,|\, \cdots\, |\, \hat{k}_{i_j}\, |\,\hat{h}_{1}\,|\, \cdots\, |\, \hat{h}_{n-i_j}\,)|.
\end{align*}
\end{claim}

\begin{proof}
We show the claim for types $\lambda=(\lambda_1,\ldots,\lambda_n)\in \inn(\sm)$ and the full standard flag $\Sc=(V_1,\ldots,V_{n})$ where $V_i=\text{span} \lbrace e_1,\ldots, e_i\rbrace$ (the $e_i$ being the standard base of $\R^n$). The claim follows then in full generality from Lemma \ref{lem gromov product continuous in types}.

Since $(\cdot\,|\,\cdot)_{I_n,\lambda}$ is invariant under the $\text{SO}(n,\R)$ action, it is enough to determine $(k\textsc{S}| \textsc{S})_{I_n,\lambda}$ or $(\textsc{S}|k \textsc{S})_{I_n,\lambda}$ for arbitrary $k\in \text{SO}(n,\R)$.
Proposition \ref{prop Gromov product as Busemann function} implies that 
$(\Sc|k \Sc)_{I_n,\lambda}=\frac{1}{2}b_{\Sc_{\lambda}}(I_n,n_{k\Sc}(I_n,\Sc)\cdot I_n)$, where $\Sc_{\lambda}$ is point in the ideal boundary $\di X(n)$ determined by the eigenvalue flag pair $(\lambda,\Sc)$ and $n_{k\Sc}(I_n,\Sc)\in N_{k\Sc}$, i.e. the element in the horospherical subgroup to $k\Sc$ such that $n_{k\Sc}(I_n,\Sc)\cdot I_n\in P(k\Sc,\Sc)$.

We first determine $n_{kS}(I_n,\Sc)\cdot I_n$. 
Let $k_w\in \mathrm{SO}(n,\R)$ be the standard antidiagonal matrix with $-1$ in the upper right corner. Then $k_w \Sc=W$ with $W$ the standard opposite flag, i.e. $W=(V^{*}_1,\ldots,V^{*}_n)$ with $V^{*}_i =\text{span} \lbrace e_n,\ldots, e_{n-i+1}\rbrace$.
Since any $k\in \text{SO}(n,\R)$ stabilizes $I_n$, the maximal flat through $k\Sc$ and $I_n$ is the unique maximal flat (i.e. affine apartment) that joins $k\Sc$ and $kW=kk_{w}S$. This yields $n_{k\Sc}(I_n,\Sc)=n_{k\Sc}(kk_{w}\Sc,\Sc)$ - here $n_{k\Sc}(kk_{w}\Sc,\Sc)\in N_{k\Sc}$ is the unique element mapping $kk_{w}\Sc$ to $\Sc$ . 

We know $N_{k\Sc}=k N_{\Sc} k^{-1}=k N_{\Sc} k^{t}$ and $N_{\Sc}$ is the group of upper triangular matrices with ones on the diagonal. Thus we are looking for $n_{\Sc}\in N_{\Sc}$ such that $k n_ {\Sc}k^{t} k k_{w} \Sc=\Sc$, i.e. $k n_\Sc k_{w}\in \st(\Sc)$; which is equivalent to $k n_{\Sc} k_{w}$ being upper triangular. 

Let $k_i$ denote the $i$-th row of $k$. Then it is straight forward to check that the $(n+1-j)$-th column of $n_s$ is given by $\sum_{i=1}^j a_{i,n+1-j} k_i$, with $a_{i,n+1-j}$ such that
\begin{align}\label{eq coefficients for bourdon fucts in type An}
\begin{pmatrix}
k_{1,n-j+1}& \cdots & & k_{j,n-j+1}\\ k_{1,n-j+2}& \cdots & & k_{j,n-j+2}\\ \vdots& \cdots & & \vdots\\ k_{1,n}& \cdots & & k_{j,n}
\end{pmatrix} 
\begin{pmatrix}
a_{1,n+1-j}\\ a_{2,n+1-j}\\ \vdots\\ a_{j,n+1-j}
\end{pmatrix}= \begin{pmatrix}
1\\ 0 \\ \vdots\\ 0
\end{pmatrix}.
\end{align}
We set $A:=k n_s$. Then $n_{k\Sc}(I_n,\Sc)\cdot I_n=(k n_{\Sc} k^t) \cdot I_n =kn_{\Sc} n_{\Sc}^t k^t= AA^t$.

The Busemann function on $X(n)$ is well known - see Lemmata 2.4, 2.5 in \cite{Hattori}. Namely, for $p\in P_n$ we have $b_{\Sc_{\lambda}}(p, I_n)=n \log (\prod_{j=1}^{n-1} (\det\Delta^{-}_{j}(p))^{\lambda_{n-j}-\lambda_{n+1-j}})$, 
where $\Delta^{-}_{j}(p)$ is the lower right $j\times j$-minor of $p$ - e.g. $\Delta^{-}_{1}(p)=p_{n,n}$. This yields $(\textsc{S}|k \textsc{S})_{I_n,\lambda}=\frac{n}{2} \sum_{j=1}^{n-1} (\lambda_{n+1-j}-\lambda_{n-j}) \log \det(\Delta^{-}_{j}(AA^t))$.

Let $(J)_{i,j}=\delta_{i,n+1-j}$ with $\delta_{i,j}$ being Kronecker's delta, i.e. $J$ is the antidiagonal. Then $AJ$ is upper triangular with the diagonal of the form $a_{1,n},\ldots, a_{n,1}$. Then one can easily show that $\Delta^{-}_{j}(AA^t)=\Delta^{-}_{j}(AJ)\Delta^{-}_{j}(JA^t)$; and thus $\det \Delta^{-}_{j}(AA^t)=\det \Delta^{-}_{j}(AJ) \det\Delta^{-}_{j}(JA^t)=a_{n,1}^2\cdots a_{n+1-j,j}^2$.

If we apply Cramer's rule to equation \eqref{eq coefficients for bourdon fucts in type An}, we get 
\begin{align*}
a_{j,n+1-j}=(-1)^{j+1}\det \begin{pmatrix}
k_{1,n-j+2}& \cdots & & k_{j-1,n-j+2}\\ \vdots& \cdots & & \vdots\\ k_{1,n}& \cdots & & k_{j-1,n}
\end{pmatrix} 
\slash \det \begin{pmatrix}
k_{1,n-j+1}& \cdots & & k_{j,n-j+1}\\ k_{1,n-j+2}& \cdots & & k_{j,n-j+2}\\ \vdots& \cdots & & \vdots\\ k_{1,n}& \cdots & & k_{j,n}
\end{pmatrix}.
\end{align*}
for $j\geq 2$ and $a_{1,n}=k_{1,n}^{-1}$. Thus $\det \Delta^{-}_{n-j}(AA^t) =
\det (\,e_1\,|\, \cdots\, |\, e_{n-j}\, |\,k_{1}\,|\, \cdots\, |\, k_{j}\,)^{2}$. 
Let $\hat{k}_i$ denote the $i$-th column of $k \in \mathrm{SO}(n,\R)$. Then
$(k\textsc{S}| \textsc{S})_{I_n,\lambda}=(\textsc{S}|k^{t} \textsc{S})_{I_n,\lambda}=n \sum_{j=1}^{n-1} (\lambda_{n+1-j}-\lambda_{n-j}) \log |\det (\,e_1\,|\, \cdots\, |\, e_{j}\, |\,\hat{k}_{1}\,|\, \cdots\, |\, \hat{k}_{n-j}\,)|,$.

Let $k,h\in \mathrm{SO}(n,\R)$. Then the $i$-th column of $h^{-1} k$ is given by $h^{-1} k\cdot e_i=h^{-1}\hat{k}_i$. Then 
$\det (\,e_1\,|\, \cdots\, |\, e_{j}\, |\,h^{-1}\hat{k}_{1}\,|\, \cdots\, |\, h^{-1}\hat{k}_{n-j}\,)=
\det (\,\hat{h}_1\,|\, \cdots\, |\, \hat{h}_{j}\, |\,\hat{k}_{1}\,|\, \cdots\, |\, \hat{k}_{n-j}\,)$
yields
$(h^{-1}k\textsc{S}|\textsc{S})_{I_n,\lambda}
=n \sum_{j=1}^{n-1} (\lambda_{n+1-j}-\lambda_{n-j}) \log|\det (\,\hat{h}_1\,|\, \cdots\, |\, \hat{h}_j\, |\,\hat{k}_{1}\,|\, \cdots\, |\, \hat{k}_{n-j}\,)|$, hence 
$(k\textsc{S}|h\textsc{S})_{I_n,\lambda}=n \sum_{j=1}^{n-1} (\lambda_{j+1}-\lambda_{j}) \log|\det (\,\hat{k}_1\,|\, \cdots\, |\, \hat{k}_j\, |\,\hat{h}_{1}\,|\, \cdots\, |\, \hat{h}_{n-j}\,)|.$
\end{proof}

\bp
Let $\lambda=(\lambda_1,\ldots,\lambda_l)$ be a type, and $\tm$ such that $\lambda\in\inn(\tm)$. Let $V=(V_1,\ldots,V_l)$, $Y=(Y_1,\ldots,Y_l)\in\Ftm$ and let $W=(W_1,\ldots,W_l)$, $Z=(Z_1,\ldots,Z_l)\in\Ftim$. Then
\begin{align*}
\cro_{\lambda}(V,W,Y,Z)= n \sum_{j=1}^{l-1}(\lambda_j - \lambda_{j+1})\log(|\frac{V_j\wedge W_{l-j}}{V_j\wedge Z_{l-j}} \frac{Y_j\wedge Z_{l-j}}{Y_j\wedge W_{l-j}}|),
\end{align*}
using the above conventions.
\ep 

\begin{proof}
As mentioned in Example \ref{ex Sl n}, the term is independent of the choices made.
By the transitivity of the $\text{SO}(n,\R)$ action, we know that every flag $V\in\Ftm$ can be written as $k\textsc{S}_{\tm}$ for $\textsc{S}_{\tm}\in\Ftm$ the standard flag and some $k\in \text{SO}(n,\R)$. Then the columns $\hat{k}_i$ are such that  $V_j=\text{span}\lbrace \hat{k}_1,\ldots, \hat{k}_{i_j} \rbrace$. In the same way every flag $W\in\Ftim$ can be written as $h\textsc{S}_{\iota\tm}$ for $\textsc{S}_{\iota\tm}\in\Ftim$ the standard flag and some $h\in \mathrm{SO}(n,\R)$.

We fix the identification $\wedge^n \R^n\simeq \det$. Then this yields $|V_j\wedge W_{l-j}|=|\det (\,\hat{k}_1\,|\, \cdots\, |\, \hat{k}_{i_j}\, |\,\hat{h}_{1}\,|\, \cdots\, |\, \hat{h}_{n-i_j}\,)|$. Thus the claim follows from the lemma above.
\end{proof}

\bibliography{mybib}{}
\bibliographystyle{abbrv}

\end{document}